\documentclass[11pt]{amsart}

\usepackage[T1]{fontenc}
\usepackage{lmodern}
\usepackage{microtype}
\usepackage[margin=1.08in,marginparwidth=0.85in,marginparsep=0.12in]{geometry}
\usepackage{amsmath,amssymb,amsthm,mathtools}
\usepackage{mathrsfs}
\usepackage{enumitem}
\usepackage{booktabs}
\usepackage{array}
\usepackage{longtable}
\usepackage{tikz-cd}
\usepackage[colorlinks=true,linkcolor=blue,citecolor=blue,urlcolor=blue]{hyperref}
\numberwithin{equation}{section}
\newtheorem{theorem}{Theorem}[section]
\newtheorem{proposition}[theorem]{Proposition}
\newtheorem{lemma}[theorem]{Lemma}
\newtheorem{corollary}[theorem]{Corollary}
\newtheorem{definition}[theorem]{Definition}

\newtheorem{maintheorem}{Theorem}

\newtheorem{maincorollary}[maintheorem]{Corollary}
\newtheorem*{restatedcorollaryB}{Corollary B}
\theoremstyle{remark}
\newtheorem{remark}[theorem]{Remark}

\title[Integral trace duality and the Witt filtration]
{Integral trace duality and the Witt filtration on Lubin--Tate space}
\author{Andrea Taccani}
\address{University of Regensburg}
\email{andrea.taccani@mathematik.uni-regensburg.de}
\date{}

\begin{document}
\raggedbottom

\begin{abstract}

Let $\mathcal O$ be the Lubin--Tate deformation ring of a one-dimensional commutative formal group of finite height $n\geq1$ over a perfect field $k$ of characteristic $p$. We identify Rezk's Witt filtration $\{I_d\}$ with the divisibility filtration of integral trace pairings on subgroup algebras of fixed source lifts. This gives the geometric recursion
\[
 I_{d+1}=pI_d+J_{d,z}\qquad(d\geq0),
\]
where $J_{d,z}$ defines the targets of order-$p^d$ quotient isogenies from a $W(k)$-valued lift $z$ of the closed point. The corrected ideal is independent of $z$.

We also compute the cokernel of the universal source trace pairing over $\mathcal O$: it is $\bigoplus_{j=1}^d(\mathcal O/p^j)^{r_j-r_{j-1}}$, where $r_j$ counts order-$p^j$ subgroups of $(\mathbf Q_p/\mathbf Z_p)^n$. The same duality extends the generic trace along two-step isogeny flags to an integral surjection, even when the flag-forgetting map is not flat. This computes composition in the subring of weighted trace operations. Applications include sharp projector denominators and congruence fibre products for full-torsion loci.
\end{abstract}

\maketitle
\tableofcontents
\section{Introduction}

Let $\Gamma/k$ be a one-dimensional commutative formal group of finite height $n\geq1$ over a perfect field of characteristic $p$. Its universal deformation space is
\[
 \mathfrak X=\operatorname{Spf}(\mathcal O),\qquad
 \mathcal O=W(k)[[u_1,\ldots,u_{n-1}]].
\]
Cofreeness equips $\mathcal O$ with Rezk's Witt filtration $\mathcal O=I_0\supset I_1\supset\cdots$. To test membership in $I_d$, one reduces coefficients to $k$ and checks that the low components of symmetric-group power operations vanish (Section~\ref{sec:2}). Rezk's intersection identity $pI_d=I_{d+1}\cap p\mathcal O$ and identification of $\mathcal O/(p,I_{d+1})$ with finite-subgroup geometry describe the intersection with $p\mathcal O$ and the reduction modulo $p$. Rezk asks for a geometric characterization of the full filtration \cite[Introduction]{RezkWitt}.

For a fixed $W(k)$-valued lift of the closed point, the integral trace pairings on its subgroup algebras determine the Witt filtration (Theorem~\ref{thm:A}). The equations of the quotient targets then give a recursion for the ideals $I_d$, in which the lift may be chosen independently at each stage (Corollary~\ref{cor:B}).

Let $\mathcal O_d$ classify subgroups of order $p^d$ of the universal deformation. Its source map $s_d:\mathcal O\to\mathcal O_d$ records the original deformation, and its target map $t_d:\mathcal O\to\mathcal O_d$ records the quotient with its induced marking. The source algebra is finite free of rank
\[
 r_d=\#\{H\subseteq(\mathbf Q_p/\mathbf Z_p)^n:|H|=p^d\},
 \qquad r_{-1}=0.
\]
The target acts on $W(k)$ by the Witt lift of the $d$th Frobenius; Section~\ref{sec:2} fixes this convention.

Fix a continuous local $W(k)$-algebra map $z:\mathcal O\to W(k)$ lifting the closed point, and put
\[
 M_{d,z}=\mathcal O_d\otimes_{s_d,\mathcal O,z}W(k),\qquad
 \tau_{d,z}(x)=t_d(x)\otimes1,\qquad
 J_{d,z}=\ker(\tau_{d,z}).
\]
Thus $J_{d,z}$ is the ideal of equations satisfied by the targets of all order-$p^d$ quotients of the fixed source deformation. Proposition~\ref{prop:3.3} shows that $\tau_{d,z}$ is surjective. In particular,
\[
 M_{d,z}\simeq\mathcal O/J_{d,z},\qquad
 J_{d,z}=\operatorname{Fitt}_0^\mathcal O(M_{d,z}),
\]
where the module action is through $\tau_{d,z}$. Write $\Theta_{d,z}(x)=\operatorname{Tr}_{M_{d,z}/W(k)}(\tau_{d,z}(x))$. After inverting the prime and splitting the subgroup algebra, this trace sums values at the quotient targets. 

\begin{maintheorem}\label{thm:A}
For every $d\geq0$, every closed-point lift $z$, and $0\leq a\leq d+1$,
\[
 I_a=\{x\in\mathcal O:
   \Theta_{d,z}(xy)\in p^aW(k)\text{ for every }y\in\mathcal O\}.
\]
If $e_{d,z}\in M_{d,z}$ is the specialized Euler class, the target map induces a ring isomorphism
\[
 \mathcal O/I_{d+1}\simeq M_{d,z}/(p e_{d,z}).
\]
\end{maintheorem}

\begin{maincorollary}\label{cor:B}
For every $d\geq0$ and every lift $z$ as above,
\begin{equation}\label{eq:1.1}
 I_{d+1}=pI_d+J_{d,z}.
\end{equation}
Choose independently a closed-point lift $z_d:\mathcal O\to W(k)$ for every $d\geq0$. Define
\begin{equation}\label{eq:1.2}
 G_0=\mathcal O,\qquad
 G_{d+1}=pG_d+\ker\!\left(
 \mathcal O\xrightarrow{t_d}\mathcal O_d\otimes_{s_d,\mathcal O,z_d}W(k)\right).
\end{equation}
Then $G_d=I_d$ for every $d$.
\end{maincorollary}

The image of $J_{d,z}$ in $\mathcal O/(pI_d)$ is $I_{d+1}/pI_d$, independently of $z$. Geometrically, the closed subscheme defined by $I_{d+1}$ is the intersection of that target locus with the subscheme defined by $pI_d$. Theorem~\ref{thm:A} identifies the rings $\mathcal O/I_{d+1}$ and $M_{d,z}/(pe_{d,z})$.

The ordinary source trace pairing in Theorem~\ref{thm:C} need not be perfect integrally. Strickland's duality, in the Euler normalization of Proposition~\ref{prop:5.2}, supplies the perfect integral functional
\[
 \lambda_d=p^{-d}\operatorname{Tr}_{(\mathcal O_d,s_d)/\mathcal O}(e_d-).
\]
Under the isomorphism $b\mapsto\lambda_d(b-)$ with the source dual, the ordinary trace map is multiplication by the complementary element $\delta_d\in\mathcal O_d$ satisfying $e_d\delta_d=p^d$ (Corollary~\ref{cor:5.4}). Its cokernel $\mathcal O_d/(\delta_d)$, annihilated by $p^d$, measures any failure of ordinary trace to supply the whole integral dual. The resulting source $\mathcal O$-module over the full deformation ring is described by Theorem~\ref{thm:C}.

\begin{maintheorem}\label{thm:C}
Let
\[
 \mathcal T_d:\mathcal O_d\longrightarrow
 \operatorname{Hom}_{\mathcal O,s_d}(\mathcal O_d,\mathcal O),\qquad
 b\longmapsto\bigl(c\mapsto\operatorname{Tr}_{(\mathcal O_d,s_d)/\mathcal O}(bc)\bigr).
\]
After independent changes of source $\mathcal O$-bases in its domain and codomain, this map is diagonal, with entry $p^j$ occurring $r_j-r_{j-1}$ times for $0\leq j\leq d$. Thus, noncanonically as an $\mathcal O$-module,
\[
 \boxed{\operatorname{coker}(\mathcal T_d)
 \simeq\bigoplus_{j=1}^{d}(\mathcal O/p^j)^{r_j-r_{j-1}}.}
\]
\end{maintheorem}

Additive Morava $E$-theory operations of weight $p^d$ are represented by source-linear functionals on $\mathcal O_d$. Theorem~\ref{thm:C} thus computes their quotient by the integral weighted traces $c\mapsto\operatorname{Tr}_{(\mathcal O_d,s_d)/\mathcal O}(bc)$, $b\in\mathcal O_d$ (Corollary~\ref{cor:8.1}). For $n>1$, $\mathcal O$ is not a principal ideal domain: the diagonal form over $\mathcal O$ does not follow from the discriminant or a single specialized Smith form.

In Section~\ref{sec:4}, transfer to the classical big Witt ring identifies target traces with ghost components, turning power-operation support into trace divisibility by $p^{d+1}$. Euler-normalized Frobenius duality identifies the elements whose trace pairing with every element of $M_{d,z}$ is divisible by $p^{d+1}$ with the ideal $pe_{d,z}M_{d,z}$. The Euler norm computes the length of $M_{d,z}/(pe_{d,z})$, matching the length given by Rezk's results, recalled in Theorem~\ref{thm:5.9} and Proposition~\ref{prop:3.1}, for $\mathcal O/I_{d+1}$. The support condition thus gives an inclusion, and surjectivity of the target map gives a surjection between these quotients of finite length. Their equal lengths force it to be an isomorphism (Proposition~\ref{prop:5.6}).

Detection of $I_{d+1}$ recovers the earlier stages through the intersection identity and gives a perfect pairing on $\mathcal O/I_{d+1}$; the trace-divisibility ideal $pe_{d,z}M_{d,z}$ gives the recursion. For the universal cokernel, ordinary trace blocks come from torsors of sections of the possibly nonsplit connected--etale extension (Section~\ref{sec:7}). Their exponents agree with the specialized Smith list, permitting the lift to a diagonal form over $\mathcal O$ by Lemma~\ref{lem:7.4}.

As an application of this duality, we describe composition of trace operations. Put $\mathcal A[d]=\operatorname{Hom}_{\mathcal O,s_d}(\mathcal O_d,\mathcal O)$ and $\mathcal A_{\mathrm{tr}}[d]=\operatorname{im}(\mathcal T_d)$. The modules $\mathcal A[d]$ are the finite-weight pieces of the additive Morava $E$-theory operation algebra. The two-step flag algebra and the flag-forgetting map are
\[
 C_{a,b}=\mathcal O_a\otimes_{t_a,\mathcal O,s_b}\mathcal O_b,
 \qquad \phi_{a,b}:\mathcal O_{a+b}\longrightarrow C_{a,b}.
\]

\begin{maintheorem}\label{thm:D}
For $a,b\geq0$, the generic relative trace of $C_{a,b}[1/p]$ over $\mathcal O_{a+b}[1/p]$ extends uniquely to an $\mathcal O_{a+b}$-linear surjection
\[
 \mathfrak t_{a,b}:C_{a,b}\longrightarrow\mathcal O_{a+b}.
\]
Its value on $1$ is a unit. For $x\in\mathcal O_a$, $y\in\mathcal O_b$, ordinary composition of operations satisfies
\[
 \mathcal T_a(x)\mathcal T_b(y)
 =\mathcal T_{a+b}\bigl(\mathfrak t_{a,b}(x\otimes y)\bigr),
 \qquad
 \mathcal A_{\mathrm{tr}}[a]\mathcal A_{\mathrm{tr}}[b]
 =\mathcal A_{\mathrm{tr}}[a+b],
\]
where the product of submodules means finite sums of products.
\end{maintheorem}

For $n\geq2$ and $a=b=1$, the flag-forgetting map fails to be flat: an order-$p^2$ cyclic group has one subgroup of order $p$, whereas an elementary abelian group of that order has $p+1$. Theorem~\ref{thm:D} constructs an integral trace across these differing generic ranks, giving a geometric formula for composition of weighted traces despite the failure of flatness. Its trace of $1$ is $w_{1,1}=1+p\varepsilon_{\mathrm{el}}$, where $\varepsilon_{\mathrm{el}}$ is the generic idempotent selecting the elementary abelian type. More generally, for $m\geq0$ and $d\geq nm$, Theorem~\ref{thm:9.2} proves that the projector onto subgroups containing $G[p^m]$ is made integral by precisely the source scalars in $p^{(n-1)m}\mathcal O$. Its fibre-product presentation specifies the congruence that functions on that locus and its complement must satisfy to extend to an integral function on the whole subgroup scheme.

\subsection*{Relation to previous work}
Rezk constructs recursive monic generators for the filtration at height two \cite{RezkWitt}. His intersection identity, characteristic-$p$ target identification and big-Witt trace formula provide the inputs for our mixed-characteristic trace characterization and the containment $J_{d,z}\subseteq I_{d+1}$. Here the recursion \eqref{eq:1.2} gives a geometric construction from the target ideal of a specified source lift, valid at arbitrary height and depth. Power operations identify the resulting ideals with the Witt filtration and prove their independence of the lift.

Strickland's subgroup identification and Frobenius duality \cite{StricklandSymmetric,StricklandDuality} supply the perfect integral pairing used here. Rezk's Hecke operations and Nelson's description of composable isogenies \cite{RezkLogarithm,NelsonWreath} provide the background for composition.

During the preparation of this paper, the author was informed that another manuscript on related questions was in circulation. At the time the present paper was submitted to arXiv, the author had not seen that manuscript, did not know the identities of its authors, and was unaware that it had been made publicly available. The results and arguments of the present paper were developed independently of that work.

Adepu, Choudhary, and Kim independently study geometric descriptions of Rezk's Witt filtration \cite[Theorem~4.2]{AdepuChoudharyKim2026}.

\subsection*{Structure of the article}
Sections~\ref{sec:2}--\ref{sec:6} prove exact detection and the geometric recursion. Section~\ref{sec:7} computes the universal trace cokernel; Sections~\ref{sec:8} and~\ref{sec:9} develop composition and type congruences. The height-two example of Section~\ref{sec:3.1} returns as a universal trace matrix in Remark~\ref{rem:7.7} and as weighted operations in Remark~\ref{rem:8.2}.

\section{Background and notation}\label{sec:2}
Fix a perfect field $k$ of characteristic $p$ and a one-dimensional commutative formal group $\Gamma/k$ of finite height $n\geq1$. Write $W=W(k)$ for the classical $p$-typical Witt ring. We recall the power-operation definition of the Witt filtration and the subgroup geometry used to detect it.

\subsection{Deformations and their classifying maps}

A test ring $B$ is a complete Noetherian local ring with maximal ideal $\mathfrak m_B$ and residue field $\kappa_B$ of characteristic $p$. A deformation to $B$ is a triple
\[
 (F,i,\alpha),\qquad i:k\longrightarrow\kappa_B,\qquad
 \alpha:F\otimes_B\kappa_B\xrightarrow{\sim}i^*\Gamma.
\]
Isomorphisms preserve the residue-field map $i$ and the marking $\alpha$. Lubin--Tate theory supplies a universal deformation $(G,\operatorname{id}_k,\alpha_G)$ \cite[Theorem~3.1]{LubinTateFormalModuli} over
\[
 \mathcal O=W[[u_1,\ldots,u_{n-1}]]
\]
and a natural bijection
\[
 \operatorname{Hom}_{\mathrm{cont,loc}}(\mathcal O,B)
 \xrightarrow{\sim}
 \left\{\begin{gathered}\text{isomorphism classes of}\\
 \text{deformations of $\Gamma$ to $B$}\end{gathered}\right\}.
\]
The marking-preserving isomorphism from a deformation classified by $f$ to $f^*G$ is unique. Since $k$ is perfect, $i$ lifts uniquely to $W\to B$.

Put $\mathfrak X=\operatorname{Spf}(\mathcal O)$, $\mathfrak m=(p,u_1,\ldots,u_{n-1})$, $\overline{\mathcal O}=\mathcal O/p\mathcal O$, and let $\epsilon:\mathcal O\twoheadrightarrow k$ be reduction modulo $\mathfrak m$.

\subsection{Restriction, transfer and total power operations}
\label{sec:2.2}
The associated Morava $E$-theory is an even-periodic, complex-orientable $E_\infty$-ring with $\pi_0E=\mathcal O$ and formal group $G$ \cite[Theorem~4.3]{RezkIsogenies}. For a finite group $H$, write $E^0BH=E^0(BH)$, with its coefficient $\mathcal O$-algebra structure. An inclusion $K\subseteq H$ gives restriction and stable transfer
\[
 \operatorname{res}_K^H:E^0BH\longrightarrow E^0BK,
 \qquad
 \operatorname{Tr}_K^H:E^0BK\longrightarrow E^0BH.
\]
Let $\Sigma_m$ be the symmetric group, with $\Sigma_0$ trivial, and use the block inclusion $\Sigma_i\times\Sigma_j\subseteq\Sigma_{i+j}$. The external product
\[
 \boxtimes:E^0B\Sigma_i\otimes_{\mathcal O}E^0B\Sigma_j
 \xrightarrow{\sim}E^0(B\Sigma_i\times B\Sigma_j)
\]
is a K\"unneth isomorphism. Finite freeness and evenness follow from \cite[Theorem~3.2]{StricklandSymmetric}; the coefficient and K\"unneth formulation for arbitrary perfect $k$ is given in \cite[Section~3.23]{RezkCongruence}. After coefficient change to a commutative $\mathcal O$-algebra $A$,
\begin{equation}\label{eq:2.2}
 (u\otimes a)\boxtimes(v\otimes b)=(u\boxtimes v)\otimes ab.
\end{equation}
Both sides take values in $E^0(B\Sigma_i\times B\Sigma_j)\otimes_{\mathcal O}A$.

The total power operation $P_m:\mathcal O\to E^0B\Sigma_m$ sends the class of $f_c:\mathbb S\to E$ to the composite

\[
 \Sigma^\infty_+B\Sigma_m =
 \Sigma^\infty_+(E\Sigma_m/\Sigma_m) \simeq
 E\Sigma_m{}_+\wedge_{\Sigma_m}\mathbb S =
 (\mathbb S^{\wedge m})_{h\Sigma_m}
 \xrightarrow{(f_c^{\wedge m})_{h\Sigma_m}}
 (E^{\wedge m})_{h\Sigma_m}
 \xrightarrow{\nu_m}E,
\]
where $h\Sigma_m$ denotes homotopy orbits and $\nu_m$ comes from the $E_\infty$ multiplication. Restriction to a point gives $c^m$, and $P_0(c)=1$, $P_1(c)=c$. The identities used below are \cite[Proof of Lemma~2.2]{RezkWitt}
\begin{equation}\label{eq:2.4}
 \begin{aligned}
 \operatorname{res}_{\Sigma_i\times\Sigma_j}^{\Sigma_{i+j}}P_{i+j}(c)
 &=P_i(c)\boxtimes P_j(c),\\
 P_m(c+c')&=\sum_{i+j=m}
 \operatorname{Tr}_{\Sigma_i\times\Sigma_j}^{\Sigma_m}
       \bigl(P_i(c)\boxtimes P_j(c')\bigr),\\
 P_m(cc')&=P_m(c)P_m(c'),\qquad P_m(1)=1.
 \end{aligned}
\end{equation}

\subsection{Compatible sequences and the Witt ring}

\begin{definition}[\texorpdfstring{$E$}{E}-typical Witt vectors]
For a commutative $\mathcal O$-algebra $A$, an element of $\mathbb W_E(A)$ is a sequence $a=(a_m)_{m\geq0}$, with $a_m\in E^0B\Sigma_m\otimes_{\mathcal O}A$, satisfying
\[
 a_0=1,\qquad
 \operatorname{res}_{\Sigma_i\times\Sigma_j}^{\Sigma_{i+j}}(a_{i+j})
 =a_i\boxtimes a_j\quad(i,j\geq0).
\]
Its ring operations are
\[
 (a+b)_m=\sum_{i+j=m}
 \operatorname{Tr}_{\Sigma_i\times\Sigma_j}^{\Sigma_m}(a_i\boxtimes b_j),
 \qquad (ab)_m=a_mb_m.
\]
The zero vector is $(1,0,0,\ldots)$ and the unit is $(1,1,1,\ldots)$.
\end{definition}

Rezk proves that these operations make $\mathbb W_E(A)$ a commutative ring \cite[Proposition~2.1]{RezkWitt}. Coefficient change acts componentwise: an $\mathcal O$-algebra map $g:A\to A'$ induces $\mathbb W_E(g)((a_m)_m)=((1\otimes g)(a_m))_m$. The identities \eqref{eq:2.4} give the ring homomorphism
\[
 P:\mathcal O\longrightarrow\mathbb W_E(\mathcal O),\qquad
 c\longmapsto(P_m(c))_{m\geq0}.
\]
Its additivity uses Witt addition. Composing with coefficient change along $\mathcal O\to A$ gives the canonical $\mathcal O$-algebra structure on $\mathbb W_E(A)$.

\begin{definition}[Witt ideals]
For $d\geq0$ put
\[
 \mathbb I_d(A)=\{a\in\mathbb W_E(A):a_m=0\text{ whenever }0<m<p^d\}.
\]
\end{definition}

These are descending ideals with $\mathbb I_0(A)=\mathbb W_E(A)$ \cite[Proposition~4.1]{RezkWitt}. Their equivalent support description is
\[
 a\in\mathbb I_d(A)
 \quad\Longleftrightarrow\quad
 a_m=0\text{ for every }m\text{ not divisible by }p^d
\]
\cite[Corollary~3.3 and Section~4]{RezkWitt}.

Equip $k$ with the $\mathcal O$-action $\epsilon$. The canonical map $\iota:\mathcal O\to\mathbb W_E(k)$ sends $c$ to $(P_m(c)\otimes1)_m$, where the tensor product uses $\epsilon$. Its degree-one component is $\epsilon(c)$; its higher components retain the reductions of the higher power operations and need not be determined by $\epsilon(c)$.

\begin{theorem}[Cofreeness, Burklund--Schlank--Yuan]
For $\Gamma$ of finite height over the perfect field $k$, the canonical structure map
\begin{equation}\label{eq:2.10}
 \iota=\mathbb W_E(\epsilon)\circ P:\mathcal O\xrightarrow{\sim}\mathbb W_E(k) 
\end{equation}
is an isomorphism of $\mathcal O$-algebras.
\end{theorem}

This is \cite[Theorem~3.4]{BSY}, in the formulation of \cite[Theorem~1.2]{RezkWitt}. Thus the full compatible sequence recovers $c$, even though its individual components have reduced coefficients. Define
\[
 I_d=\iota^{-1}(\mathbb I_d(k)),\qquad
 \overline I_d=(I_d+p\mathcal O)/p\mathcal O.
\]
In particular, $c\in I_d$ means that $P_m(c)\otimes1=0$ for $0<m<p^d$. The reduction of $I_d$ has kernel $I_d\cap p\mathcal O$. The classical $W(k)$ is the coefficient subring of $\mathcal O$, whereas $\mathbb W_E(k)$ is the deformation ring itself. At height one, $I_d=p^d\mathcal O$. In higher height, reduction modulo $p$ gives a nontrivial filtration in the Lubin--Tate deformation directions.

\subsection{The subgroup algebra and the source map}

For $d\geq0$, Strickland's subgroup theorem \cite[Theorem~10.1]{StricklandSubgroups} gives a complete local ring $\mathcal O_d$ classifying pairs $((F,i,\alpha),H)$, where $H\subseteq F$ is finite locally free of rank $p^d$. The forgetful map defines the source structure $s_d:\mathcal O\to\mathcal O_d$, for which $\mathcal O_d$ is finite free \cite[Theorem~10.1]{StricklandSubgroups}. For a classifying map $f:\mathcal O\to B$, representability gives
\begin{equation}\label{eq:2.12}
 \begin{gathered}
 \{h:\mathcal O_d\to B\text{ continuous local}:h\circ s_d=f\}\\
 \updownarrow\ {\scriptstyle\sim}\\[-0.3ex]
 \{H\subseteq f^*G:\operatorname{rank}_B(H)=p^d\}.
 \end{gathered}
\end{equation}
The identity of $\mathcal O_d$ classifies the universal subgroup $H_d\subseteq s_d^*G$; a point $h$ classifies its pullback $h^*H_d$.

\subsection{The quotient and the target map}

A finite locally free subgroup $H\subseteq F$ of rank $p^d$ has a quotient formal group and quotient isogeny $q:F\to F/H$ \cite[Theorem~5.3]{StricklandSubgroups}, \cite[Proposition~11.9]{RezkCongruence}. Write $\phi:k\to k$, $a\mapsto a^p$, put $\Gamma^{(d)}=(\phi^d)^*\Gamma$, and let $\varphi^d:\Gamma\to\Gamma^{(d)}$ be relative Frobenius. Its coordinate pullback sends $X$ to $X^{p^d}$.

The special fibre of $H$ has ideal $(X^{p^d})$ in the discrete valuation ring $\kappa_B[[X]]$, since its colength is $p^d$. It is therefore the kernel of relative Frobenius. The quotient property determines the isomorphism $\beta$ in
\[
 \begin{tikzcd}[column sep=large,row sep=large]
 F_{\kappa_B} \arrow[r,"\overline q"] \arrow[d,"\alpha"',"\sim"]
   & (F/H)_{\kappa_B} \arrow[d,"\beta","\sim"'] \\
 i^*\Gamma \arrow[r,"i^*\varphi^d"']
   & i^*\Gamma^{(d)}=(i\circ\phi^d)^*\Gamma .
 \end{tikzcd}
\]
Thus the target is the deformation
\begin{equation}\label{eq:2.14}
 (F/H,\ i\circ\phi^d,\ \beta).
\end{equation}
This construction is natural in $B$ and gives the target map in
\[
 s_d,t_d:\mathcal O\longrightarrow\mathcal O_d.
\]
At a point $h:\mathcal O_d\to B$, $h\circ s_d$ classifies the source and $h\circ t_d$ the marked quotient. The universal quotient is
\[
 q_d:s_d^*G\longrightarrow t_d^*G,\qquad \ker(q_d)=H_d.
\]

\subsection{The cohomological description of the two maps}
\label{sec:2.6}

Let $\mathcal I_d^{\mathrm{tr}}\subseteq E^0B\Sigma_{p^d}$ be the ideal generated by the transfers from $\Sigma_i\times\Sigma_{p^d-i}$ for $0<i<p^d$. Strickland's subgroup theorem identifies
\begin{equation}\label{eq:2.17}
 \mathcal O_d=E^0B\Sigma_{p^d}/\mathcal I_d^{\mathrm{tr}},
 \qquad \operatorname{rank}_{\mathcal O,s_d}(\mathcal O_d)=N_n(p^d),
\end{equation}
where $N_n(p^d)$ counts order-$p^d$ subgroups of $(\mathbf Q_p/\mathbf Z_p)^n$ \cite[Theorem~1.1]{StricklandSymmetric}. 

Write $\kappa_d^{\mathrm{tr}}:E^0B\Sigma_{p^d}\twoheadrightarrow\mathcal O_d$ for the quotient and $\eta_{p^d}:\mathcal O\to E^0B\Sigma_{p^d}$ for the coefficient map. Then
\begin{equation}\label{eq:2.18}
 s_d=\kappa_d^{\mathrm{tr}}\circ\eta_{p^d},\qquad
 t_d=\kappa_d^{\mathrm{tr}}\circ P_{p^d}.
\end{equation}
Indeed, the mixed terms in the addition formula \eqref{eq:2.4} lie in the transfer ideal, so $P_{p^d}$ becomes a ring map in the quotient. Its identification with the framed geometric target is due to Ando--Hopkins--Strickland \cite[Proposition~12.12]{AHSSigma}. 

\begin{remark}[The coefficient action of the target map]
\label{rem:2.4}
Put $\phi^d:k\to k$ and $\sigma_d=W(\phi^d):W\to W$. By \eqref{eq:2.14} and uniqueness of the coefficient lift,
\[
 t_d(c)=s_d(\sigma_d(c))\qquad(c\in W).
\]
Consequently $t_d(cx)=s_d(\sigma_d(c))t_d(x)$: the target map is $\sigma_d$-semilinear for the source coefficient structure. Perfectness makes $\sigma_d$ an automorphism, and $\sigma_d(p)=p$.

Traces, determinants, duals and source base change use the source structure unless explicitly stated otherwise. A source-$W$-linear functional $\ell$ satisfies $\ell(t_d(cx))=\sigma_d(c)\ell(t_d(x))$; applying $\sigma_d^{-1}$ to its values restores $W$-linearity in $x$. This transport preserves $p$-divisibility and module lengths; compare \cite[Corollary~3.49]{BSY}. Here $\sigma_d$ always denotes the coefficient automorphism, distinct from BSY's source projection and from the Witt transfer $\sigma$ of Section~\ref{sec:4}.
\end{remark}

\section{Target loci of fixed source lifts}\label{sec:3}
Fix a continuous local $W$-algebra map $z:\mathcal O\to W$ reducing to the closed point, or equivalently choose $z(u_i)=c_i\in pW$. Let $G_z=z^*G$. Fixing this source while varying its subgroup gives the finite free $W$-algebra
\[
 M_{d,z}=\mathcal O_d\otimes_{s_d,\mathcal O,z}W,
 \qquad \operatorname{rank}_W M_{d,z}=N_n(p^d),
\]
with relations
\[
 s_d(a)\otimes1=1\otimes z(a)\qquad(a\in\mathcal O).
\]
For a complete local $W$-algebra $B$, continuous $W$-algebra maps $M_{d,z}\to B$ correspond by \eqref{eq:2.12} to rank-$p^d$ finite locally free subgroups of $(G_z)_B$. Write $\pi_{d,z}:\mathcal O_d\to M_{d,z}$ for $b\mapsto b\otimes1$ and define
\[
 \tau_{d,z}=\pi_{d,z}\circ t_d:\mathcal O\longrightarrow M_{d,z},
 \qquad \tau_{d,z}(x)=t_d(x)\otimes1.
\]
For a subgroup classified by $h:M_{d,z}\to B$, the composite $h\circ\tau_{d,z}$ classifies its marked quotient. Hence
\[
 J_{d,z}=\ker(\tau_{d,z})
\]
is the ideal of equations satisfied by these targets. Through the source coefficient action $c\mapsto1\otimes c$, one has
\[
 \tau_{d,z}(cx)=\sigma_d(c)\tau_{d,z}(x)
 \qquad(c\in W,\ x\in\mathcal O).
\]

The closed source fibre is $R_d=\mathcal O_d\otimes_{s_d,\mathcal O,\epsilon}k$. Since $z\bmod p=\epsilon$, base change gives $M_{d,z}/pM_{d,z}\simeq R_d$, independently of $z$. Reduction of the target gives
\[
 \begin{tikzcd}[column sep=large,row sep=large]
 \mathcal O \arrow[r,"\tau_{d,z}"] \arrow[d,two heads]
   & M_{d,z} \arrow[d,two heads] \\
 \overline{\mathcal O} \arrow[r,"\overline t_d"'] & R_d .
 \end{tikzcd}
\]
The horizontal maps are semilinear for the source coefficient actions, with automorphisms $\sigma_d$ and $\phi^d$, respectively. 

\begin{proposition}\label{prop:3.1}
For every $d\geq0$, reduction of $t_d$ induces a surjective ring map
\[
 \overline t_d:\overline{\mathcal O}\twoheadrightarrow R_d,
 \qquad \ker(\overline t_d)=\overline I_{d+1}.
\]
\end{proposition}

\begin{proof}
Surjectivity is \cite[Proposition~8.1]{RezkWitt}. The $d$th ghost at $k$ is
\[
 g_d:\mathbb W_E(k)\longrightarrow R_d,\qquad
 (a_m)_m\longmapsto(\kappa_d^{\mathrm{tr}}\otimes1)(a_{p^d}),
\]
where the tensor product uses $\epsilon$. Equations \eqref{eq:2.10} and \eqref{eq:2.18} give $g_d(\iota(c))=t_d(c)\otimes1=\overline t_d(c\bmod p)$. Applied to the perfect $k$-algebra $k$, \cite[Proposition~8.8]{RezkWitt} gives $\ker(g_d)=p\mathbb W_E(k)+\mathbb I_{d+1}(k)$. Transport through $\iota$ and then reduction of the domain modulo $p$ yield
\[
 \ker(\overline t_d)=(p\mathcal O+I_{d+1})/p\mathcal O=\overline I_{d+1}.
\]
\end{proof}

Thus the characteristic-$p$ target coordinates present the subgroup fibre as $\overline{\mathcal O}/\overline I_{d+1}\simeq R_d$. We recall the elementary lifting argument needed in mixed characteristic.

\begin{lemma}\label{lem:3.2}
Let $f:B\to C$ be a continuous homomorphism of $p$-adically complete and separated rings. If $B/pB\to C/pC$ is surjective, then $f$ is surjective.
\end{lemma}

\begin{proof}
For $c\in C$, choose recursively $b_r\in B$ so that
\[
 c-f\!\left(\sum_{i<r}p^ib_i\right)\in p^rC\quad(r\geq0).
\]
To pass from $r$ to $r+1$, divide this error by $p^r$ in $C$ and lift its reduction modulo $p$ to $B$. The sum $b=\sum_{r\geq0}p^rb_r$ converges in $B$; continuity and separation of $C$ give $f(b)=c$.
\end{proof}

\begin{proposition}\label{prop:3.3}
For every $d\geq0$ and every continuous local $W$-algebra map $z:\mathcal O\to W$ reducing to the closed point, $\tau_{d,z}$ is surjective and
\begin{equation}\label{eq:3.6}
 M_{d,z}\simeq\mathcal O/J_{d,z},\qquad
 J_{d,z}=\operatorname{Fitt}^{\mathcal O}_0(M_{d,z}),
\end{equation}
where the $\mathcal O$-module structure is the target action through $\tau_{d,z}$. Moreover,
\begin{equation}\label{eq:3.7}
 J_{d,z}/pJ_{d,z}\xrightarrow{\sim}\overline I_{d+1},
\end{equation}
and
\begin{equation}\label{eq:3.8}
 J_{d,z}\cap p\mathcal O=pJ_{d,z}.
\end{equation}
\end{proposition}

\begin{proof}
Both $\mathcal O=W[[u_1,\ldots,u_{n-1}]]$ and the finite free $W$-module $M_{d,z}$ are $p$-adically complete and separated. Every ring map preserves powers of $p$, so $\tau_{d,z}$ is also continuous for these topologies. Proposition~\ref{prop:3.1} and Lemma~\ref{lem:3.2} give its surjectivity. If $j_1,\ldots,j_r$ generate $J_{d,z}$, the presentation $\mathcal O^r\xrightarrow{(j_1\ \cdots\ j_r)}\mathcal O\to M_{d,z}\to0$ gives the Fitting equality in \eqref{eq:3.6}.

If $px\in J_{d,z}$, then $p\tau_{d,z}(x)=0$. The $p$-torsion-freeness of $M_{d,z}$ forces $x\in J_{d,z}$, proving \eqref{eq:3.8}. Reducing the quotient presentation and applying Proposition~\ref{prop:3.1} now gives
\[
 J_{d,z}/pJ_{d,z}\simeq(J_{d,z}+p\mathcal O)/p\mathcal O
 =\overline I_{d+1},
\]
which proves \eqref{eq:3.7}.
\end{proof}

\subsection{Running example at height two}
\label{sec:3.1}

Take $p=2$, $k=\mathbf F_2$, and the height-two formal group of the supersingular elliptic curve $y^2+y=x^3$. In the model of \cite[Sections~3--4]{RezkHeightTwo},
\[
 \mathcal O=\mathbf Z_2[[u]],\qquad
 \mathcal O_1=\mathcal O[b]/(b^3-ub-2),
\]
and
\begin{equation}\label{eq:3.9}
 s_1(u)=u,\qquad t_1(u)=u^2+3b-ub^2.
\end{equation}
The coefficient Frobenius is trivial. For $z_0(u)=0$,
\[
 M_{1,z_0}=\mathbf Z_2[b]/(b^3-2),\qquad
 \tau_{1,z_0}(u)=3b.
\]
The basis $1,3b,9b^2$ identifies $\mathbf Z_2[[u]]/(u^3-54)$ with $M_{1,z_0}$, so
\[
 J_{1,z_0}=(u^3-54).
\]
Using $I_1=(2,u)$, Theorem~\ref{thm:6.1} gives
\[
 I_2=2I_1+J_{1,z_0}=(4,2u,u^3-2).
\]
The correction $2I_1$ absorbs the difference between the target equation $u^3-54$ and the displayed Witt generator. We return to this example in Remark~\ref{rem:7.7}.

\section{Trace divisibility from the Witt support}\label{sec:4}
Fix a closed-point lift $z:\mathcal O\to W(k)$ as in Section~\ref{sec:3}. For $r\geq0$, define the target trace
\[
 \Theta_{r,z}(x)
 :=\operatorname{Tr}_{M_{r,z}/W(k)}\bigl(m_{\tau_{r,z}(x)}\bigr)
 \qquad(x\in\mathcal O),
\]
where $m_b$ is multiplication by $b$ on the finite free source $W(k)$-module $M_{r,z}$. Remark~\ref{rem:2.4} gives
\[
 \Theta_{r,z}(cx)=\sigma_r(c)\Theta_{r,z}(x)
 \qquad(c\in W(k),\ x\in\mathcal O).
\]

For a commutative $\mathcal O$-algebra $A$ and $d\geq0$, the $d$th $E$-typical ghost \cite[Section~5]{RezkWitt} takes one component and then passes to the transfer quotient. It is the ring map
\begin{equation}\label{eq:4.2}
 \operatorname{gh}_d:\mathbb W_E(A)\longrightarrow
 \mathcal O_d\otimes_{s_d,\mathcal O}A,
 \qquad (a_m)_m\longmapsto a_{p^d}\bmod \mathcal I_d^{\mathrm{tr}}.
\end{equation}
It is generally not $\mathcal O$-linear for the source action. The target construction in Section~\ref{sec:2.6} gives
\begin{equation}\label{eq:4.3}
 t_d=\operatorname{gh}_d\circ P:\mathcal O\longrightarrow\mathcal O_d;
\end{equation}
see \cite[Section~8, before Proposition~8.1]{RezkWitt}.

Write $\mathbb W_{\mathbf Z}(A)$ for the classical big Witt ring, with power-series coordinate
\[
 \mathscr S:(\mathbb W_{\mathbf Z}(A),+)\xrightarrow{\sim}(1+XA[[X]],\times).
\]
For $F(X)=\mathscr S(w)$ and $N\geq1$, its ghost is
\[
 \operatorname{gh}_N(w)=[X^N]\frac{XF'(X)}{F(X)}.
\]
Here $[X^N]$ denotes coefficient extraction. The logarithmic derivative turns multiplication of power series into addition of ghost components; it is defined over every $A$, since $F(0)=1$. Over $W(k)$, we have
\begin{equation}\label{eq:4.4}
 \operatorname{gh}_N(w)=N[X^N]\log\mathscr S(w),
\end{equation}
with the formal logarithm computed in $W(k)[1/p][[X]]$; see \cite[Section~9.1]{RezkWitt}.

For a finite group $Q$, put $LQ=L_{K(n)}\Sigma^\infty_+BQ$ and let $e$ be the trivial group. Transposing $LQ\to Le$ under Strickland's canonical self-duality \cite[Proposition~8.3 and Definition~8.4]{StricklandDuality} and applying $E^0$ defines the $\mathcal O$-linear local transfer
\begin{equation}\label{eq:4.5}
 \epsilon_Q=(Q\to e)_!:E^0BQ\longrightarrow\mathcal O.
\end{equation}
For $Q=\Sigma_m$, write $\operatorname{Tr}^{e}_{\Sigma_m}=\epsilon_{\Sigma_m}$, and extend this map $A$-linearly after coefficient change. This local transfer should be distinguished from augmentation and the ordinary algebra trace.

\begin{proposition}
For every commutative $\mathcal O$-algebra $A$, there is a natural additive homomorphism
\[
 \sigma_A:\mathbb W_E(A)\longrightarrow\mathbb W_{\mathbf Z}(A)
\]
such that, for $a=(a_m)_m\in\mathbb W_E(A)$,
\begin{equation}\label{eq:4.6}
 \mathscr S(\sigma_A(a))
 =\sum_{m\geq0}\operatorname{Tr}^{e}_{\Sigma_m}(a_m)X^m,
\end{equation}
and, for $r\geq0$,
\begin{equation}\label{eq:4.7}
 \operatorname{gh}_{p^r}(\sigma_A(a))
 =\operatorname{Tr}_{(\mathcal O_r\otimes_{s_r,\mathcal O}A)/A}
       \bigl(m_{\operatorname{gh}_r(a)}\bigr).
\end{equation}
\end{proposition}

These formulas are \cite[Propositions~11.2 and~11.1]{RezkWitt}, respectively. The algebra in \eqref{eq:4.7} is finite free over $A$ by source base change. Take $A=W(k)$ via $z$, write $\sigma=\sigma_{W(k)}$, and apply the proposition to $a=\mathbb W_E(z)(P(x))$. Naturality of \eqref{eq:4.2} and \eqref{eq:4.3} identifies its $r$th ghost with $t_r(x)\otimes1=\tau_{r,z}(x)$, giving
\begin{equation}\label{eq:4.8}
 \operatorname{gh}_{p^r}\!\left(\sigma\bigl(\mathbb W_E(z)(P(x))\bigr)\right)
 =\operatorname{Tr}_{M_{r,z}/W(k)}\bigl(m_{\tau_{r,z}(x)}\bigr)
 =\Theta_{r,z}(x).
\end{equation}

Lemmas~\ref{lem:4.2} and~\ref{lem:4.3} supply the valuation estimate and the componentwise divisibility needed to apply it.

\begin{lemma}\label{lem:4.2}
Let $a\geq1$, $r\geq a-1$, and $N=p^r$. If
\[
 F(X)=1+\sum_{m\geq1}b_mX^m\in1+XW(k)[[X]]
\]
satisfies $b_m\in pW(k)$ for $1\leq m<p^a$, then
\[
 N[X^N]\log F(X)\in p^aW(k).
\]
\end{lemma}

\begin{proof}
Put $B=F-1$. Since $B(0)=0$, only $1\leq j\leq N$ contributes to $[X^N]$ in
\[
 N\log(1+B)=\sum_{j\geq1}(-1)^{j+1}\frac NjB^j.
\]
Consider a term
\[
 \frac Nj b_{m_1}\cdots b_{m_j},
 \qquad m_i\geq1,\quad m_1+\cdots+m_j=N.
\]
Let $\ell=\#\{i:m_i<p^a\}$ and $v=v_p(j)$. The product contributes at least $\ell$ factors of $p$, while $N/j$ has valuation $r-v$. Thus the term has valuation at least $r-v+\ell$.

If $\ell=j$, then $j-v_p(j)\geq1$, giving $r-v+\ell\geq r+1\geq a$. If $\ell<j$, some $m_i\geq p^a$, hence $r\geq a$, and
\[
 p^r\geq(j-\ell)p^a+\ell,
 \qquad j\leq p^{r-a}+\ell.
\]
Thus $v\leq r-a+\ell$: otherwise
\[
 j\geq p^{r-a+\ell+1}>p^{r-a}+\ell,
\]
a contradiction. The strict inequality uses $p\geq2$ and $2^{\ell+1}\geq\ell+2$. Again $r-v+\ell\geq a$, so every contribution lies in $p^aW(k)$.
\end{proof}

\begin{lemma}\label{lem:4.3}
For $y\in\mathcal O$, put
\[
 \eta_y=\mathbb W_E(z)(P(y))=(\eta_{y,m})_{m\geq0}\in\mathbb W_E(W(k)).
\]
If $a\geq0$ and $y\in I_a$, then
\[
 \eta_{y,m}\in p\bigl(E^0B\Sigma_m\otimes_{\mathcal O,z}W(k)\bigr)
 \qquad(0<m<p^a).
\]
\end{lemma}

\begin{proof}
Let $\pi:W(k)\to k$ be reduction. Since $\pi\circ z=\epsilon$, functoriality and \eqref{eq:2.10} give
\[
 \mathbb W_E(\pi)(\eta_y)=\mathbb W_E(\epsilon)(P(y))=\iota(y).
\]
On the $m$th coefficient module $N_m:=E^0B\Sigma_m\otimes_{\mathcal O,z}W(k)$, this is
\[
 N_m\longrightarrow E^0B\Sigma_m\otimes_{\mathcal O,\epsilon}k
 \simeq N_m/pN_m.
\]
For $y\in I_a$, the components of $\iota(y)$ with $0<m<p^a$ vanish by definition of the Witt filtration. Hence $\eta_{y,m}\in pN_m$.
\end{proof}

\begin{lemma}\label{lem:4.4}
Let $a\geq1$ and $r\geq a-1$. Then
\[
 \boxed{\Theta_{r,z}(I_a)\subseteq p^aW(k).}
\]
\end{lemma}

\begin{proof}
For $y\in I_a$, use $\eta_y$ from Lemma~\ref{lem:4.3} and write
\[
 F_y(X)=\mathscr S(\sigma(\eta_y))
       =1+\sum_{m\geq1}b_mX^m,
 \qquad b_m=\operatorname{Tr}^{e}_{\Sigma_m}(\eta_{y,m}),
\]
using \eqref{eq:4.6}. Lemma~\ref{lem:4.3} and $W(k)$-linearity of the transfer give $b_m\in pW(k)$ for $0<m<p^a$. Lemma~\ref{lem:4.2}, together with \eqref{eq:4.4} and \eqref{eq:4.8}, therefore yields
\[
 \Theta_{r,z}(y)
 =\operatorname{gh}_{p^r}(\sigma(\eta_y))
 =p^r[X^{p^r}]\log F_y(X)\in p^aW(k).
\]
\end{proof}

\begin{remark}

The range $r\geq a-1$ in Lemma~\ref{lem:4.4} cannot be weakened in general. In the height-two model at $p=2$ of Section~\ref{sec:3.1}, with $z_0(u)=0$, one has
\[
 I_2=(4,2u,u^3-2),\qquad
 M_{1,z_0}=\mathbf Z_2[b]/(b^3-2),\qquad
 \tau_{1,z_0}(u)=3b;
\]
the formula for $I_2$ is also \cite[Example~8.14]{RezkWitt}. Since $M_{0,z_0}=\mathbf Z_2$ and $\tau_{0,z_0}=z_0$, the element $y=u^3-2\in I_2$ satisfies
\[
 \Theta_{0,z_0}(y)=-2\notin4\mathbf Z_2.
\]
Nor does divisibility of a single scalar trace imply membership in $I_a$: in the same example,
\[
 \Theta_{1,z_0}(u)=\operatorname{Tr}(m_{3b})=0,
 \qquad u\notin I_2,
\]
because the reduction of $I_2$ modulo $2$ is $(u^3)$. Multiplication by $u^2$ exposes the failure:
$\Theta_{1,z_0}(u\cdot u^2)=\operatorname{Tr}(m_{27b^3})=162\notin4\mathbf Z_2$.
Thus Lemma~\ref{lem:4.4} gives the stated forward implication; a converse requires conditions on traces of products.
\end{remark}

\section{Trace duality and exact detection of the Witt filtration}\label{sec:5}
We use integral duality on the subgroup algebras to characterize the Witt filtration by divisibility of traces of products.
Unless explicitly stated otherwise, every $\mathcal O$-module structure, dual, norm and trace involving $\mathcal O_d$ uses the source map $s_d:\mathcal O\to \mathcal O_d$.
For $d\geq0$ put
\[
  r_d=\operatorname{rank}_{\mathcal O}(\mathcal O_d)=N_n(p^d),
 \qquad r_{-1}=0.
\]
Counting finite-index lattices in $\mathbf Z_p^n$ gives the generating series \cite[Corollary~10.3]{StricklandSubgroups}
\begin{equation}\label{eq:5.1}
 \sum_{d\geq0}r_dT^d
 =\prod_{i=0}^{n-1}(1-p^iT)^{-1}.
\end{equation}

\subsection{Perfect duality}

Fix a degree-zero complex orientation. Put
\[
 R_d^{\Sigma}=E^0B\Sigma_{p^d},\qquad
 c_d^{\Sigma}=e(V_{p^d}-1)\in R_d^{\Sigma},
\]
where $V_{p^d}$ is the complex permutation representation. Retain the partition-transfer ideal $\mathcal I_d^{\mathrm{tr}}$ of \eqref{eq:2.17}, and let $\mathcal J_d^{\Sigma}$ be the intersection of the kernels of restriction to $\Sigma_i\times\Sigma_{p^d-i}$ for $0<i<p^d$. Strickland proves \cite[Theorem~8.6]{StricklandSymmetric}
\begin{equation}\label{eq:5.2}
 \begin{gathered}
 \mathcal O_d=R_d^{\Sigma}/\mathcal I_d^{\mathrm{tr}},\qquad
 \mathcal I_d^{\mathrm{tr}}
   =\operatorname{Ann}(c_d^{\Sigma})
   =\operatorname{Ann}(\mathcal J_d^{\Sigma}),\\
 \mathcal J_d^{\Sigma}=c_d^{\Sigma}R_d^{\Sigma}
          \simeq \mathcal O_d.
 \end{gathered}
\end{equation}
Both ideals are $\mathcal O$-module direct summands of $R_d^{\Sigma}$ \cite[Theorem~8.6]{StricklandSymmetric}; the identification with the subgroup ring is \cite[Theorem~9.2]{StricklandSymmetric}. The isomorphism $\mathcal O_d\to\mathcal J_d^\Sigma$ in \eqref{eq:5.2} is multiplication by $c_d^\Sigma$. Write $e_d$ for the image of $c_d^{\Sigma}$ in $\mathcal O_d$.
The coefficient extension that applies these statements to our arbitrary perfect field $k$ is explained in Remark~\ref{rem:5.3}.

For a finite projective $A$-algebra $R$, a \emph{Frobenius functional} is an $A$-linear map $\ell:R\to A$ such that
\[
 R\longrightarrow\operatorname{Hom}_A(R,A),\qquad
 x\longmapsto(y\mapsto\ell(xy))
\]
is an isomorphism; every $A$-linear functional on $R$ can be written uniquely as $y\mapsto\ell(xy)$. Lemma~\ref{lem:5.1} explains how a quotient of $R$ inherits such a functional when its annihilator has a suitable generator.

\begin{lemma}
\label{lem:5.1}
Let $R$ be a finite projective commutative algebra over a ring $A$ and let $\ell:R\to A$ have perfect pairing $(x,y)\mapsto\ell(xy)$. Let $I\subset R$ be an ideal such that $R/I$ is finite projective over $A$. Then
\begin{equation}\label{eq:5.3}
 \operatorname{Ann}_R(I)\xrightarrow{\sim}\operatorname{Hom}_A(R/I,A),
 \qquad a\longmapsto(\overline b\mapsto\ell(ab)).
\end{equation}
If multiplication by $c\in R$ induces an isomorphism $R/I\xrightarrow{\sim}\operatorname{Ann}_R(I)$, the functional $\overline b\mapsto\ell(cb)$ is perfect on $R/I$.
\end{lemma}

\begin{proof}
The displayed functional is independent of the lift $b$ because $aI=0$. Conversely, precompose a functional on $R/I$ with $R\to R/I$ and use perfectness of $\ell$ to represent it by a unique $a\in R$. For $i\in I$ and $r\in R$, its value on $ir$ is zero, so $\ell((ai)r)=0$ for every $r$. Perfectness implies $ai=0$. This proves \eqref{eq:5.3}; composing it with multiplication by $c$ proves the last assertion.
\end{proof}

\begin{proposition}\label{prop:5.2}
The formula
\begin{equation}\label{eq:5.4}
 \lambda_d(b):=p^{-d}\operatorname{Tr}_{\mathcal O_d/\mathcal O}(e_db)
 \qquad(b\in \mathcal O_d)
\end{equation}
defines an $\mathcal O$-linear map $\mathcal O_d\to\mathcal O$, and
\begin{equation}\label{eq:5.5}
 \mathcal O_d\xrightarrow{\sim}\operatorname{Hom}_{\mathcal O,s_d}(\mathcal O_d,\mathcal O),
 \qquad b\longmapsto(c\mapsto\lambda_d(bc)).
\end{equation}
\end{proposition}

\begin{proof}
Let $\epsilon_{\Sigma}=\epsilon_{\Sigma_{p^d}}:R_d^\Sigma\to\mathcal O$ be the local transfer of \eqref{eq:4.5}, so $\epsilon_\Sigma=\operatorname{Tr}^{e}_{\Sigma_{p^d}}$. Under Strickland's Frobenius-object structure this is the Frobenius counit, and its pairing is perfect \cite[Theorem~8.7 and the discussion following Remark~9.1]{StricklandDuality}, using the finite free even cohomology recalled in Section~\ref{sec:2.2}. For a Frobenius functional $\ell:R\to A$ and an ideal $K\subset R$, write $K^\perp=\{r\in R:\ell(rk)=0\text{ for every }k\in K\}$. Then $K^\perp=\operatorname{Ann}(K)$; since $\mathcal J_d^{\Sigma}$ is a module direct summand, $(\mathcal J_d^{\Sigma})^{\perp\perp}=\mathcal J_d^{\Sigma}$. Thus
\[
 \operatorname{Ann}(\mathcal I_d^{\mathrm{tr}})
 =\mathcal J_d^{\Sigma}.
\]
Together with \eqref{eq:5.2} and Lemma~\ref{lem:5.1}, this proves that
\[
 b\longmapsto\epsilon_{\Sigma}(c_d^{\Sigma}\widetilde b)
\]
is a well-defined integral perfect functional on $\mathcal O_d$; here $\widetilde b$ is any lift to $R_d^{\Sigma}$.

Let $C=D'$ be the faithfully flat $\mathcal O[1/p]$-algebra constructed immediately before \cite[Theorem~10.2]{StricklandDuality}; the theorem, compatibly with transfers, gives
\[
 C\otimes_{\mathcal O}R_d^\Sigma\simeq\prod_{[S]}C,
\]
where $[S]$ runs through isomorphism classes of finite $\mathbf Z_p^n$-sets of cardinality $p^d$. 

The character of a transfer from $\Sigma_i\times\Sigma_{p^d-i}$ is a sum over invariant subsets of cardinality $i$. It is therefore zero at a transitive action. Conversely, a nontransitive action has a nonempty proper invariant subset. For its cardinality $i$, the transfer of $1$ has a positive integer character value, equal to the number of such subsets. This integer is a unit in $C$. Multiplying by the idempotent of that coordinate shows that the transfer ideal contains that factor. Therefore the quotient consists of the transitive factors.

Put $\Lambda=(\mathbf Q_p/\mathbf Z_p)^n$. A transitive action is $\mathbf Z_p^n/L$ for an index-$p^d$ lattice $L$; its Pontryagin dual is an order-$p^d$ subgroup $H\subset\Lambda$. Hence
\begin{equation}\label{eq:5.7}
 C\otimes_{\mathcal O}\mathcal O_d
 \simeq\prod_{\substack{H\subset\Lambda\\|H|=p^d}}C.
\end{equation}
We write $b(H)$ and $e_d(H)$ for the corresponding coordinates. Notice that $\mathcal O_d[1/p]$ is finite etale: this can be checked after the faithfully flat extension to $C$, where it is a product of copies of $C$.

The Frobenius character formula preceding \cite[Theorem~10.2]{StricklandDuality} weights the value at $S$ by $1/|\operatorname{Aut}_{\mathbf Z_p^n}(S)|$. The character of $c_d^\Sigma$ is zero on every nontransitive $S$, because the reduced permutation representation still has a nonzero trivial summand. For a transitive $S$, the abelian image of $\mathbf Z_p^n$ acts regularly. An equivariant permutation is a translation, determined by the image of one point, so there are exactly $p^d$ such permutations. Thus, over $C$,
\[
 \epsilon_{\Sigma}(c_d^{\Sigma}\widetilde b)
   =\sum_{\substack{H\subset(\mathbf Q_p/\mathbf Z_p)^n\\|H|=p^d}}
      \frac{e_d(H)b(H)}{p^d}
   =p^{-d}\operatorname{Tr}_{\mathcal O_d/\mathcal O}(e_db).
\]
The last equality uses \eqref{eq:5.7}: the trace of multiplication in a product algebra is the sum of its coordinates. Faithful flatness gives the same equality over $\mathcal O[1/p]$. The left side is integral by construction, and $\mathcal O\to\mathcal O[1/p]$ is injective, thus the rational expression identifies the integral functional itself.
\end{proof}

\begin{remark}\label{rem:5.3}
We explain the coefficient extension from the conventions of \cite{StricklandSymmetric} to an arbitrary perfect field $k$. Choose an algebraic closure $\bar k$ and compatible deformation coordinates:
\[
 \mathcal O=W(k)[[u_1,\ldots,u_{n-1}]]
 \longrightarrow
 \mathcal O_{\bar k}=W(\bar k)[[u_1,\ldots,u_{n-1}]].
\]
Both rings are Noetherian and local. The Koszul resolution on $p,u_1,\ldots,u_{n-1}$ remains exact over $\mathcal O_{\bar k}$, so $\operatorname{Tor}^{\mathcal O}_1(k,\mathcal O_{\bar k})=0$. Applying the local flatness criterion to $\mathcal O\to\mathcal O_{\bar k}$ with the target as its own finite module proves flatness \cite[Lemma~10.99.7, Tag~00MK]{StacksFlatness}. Since the map is local, it is faithfully flat.

The functorial Goerss--Hopkins--Miller construction \cite[Theorem~4.2]{BehrensConstruction} realizes this coefficient map by $E\to E_{\bar k}$; see also \cite[Definition~2.32 and Theorem~2.38(3)]{BSY} for this perfect-field base change. Rezk's finite free coefficient formula therefore gives
\[
 \mathcal O_{\bar k}\otimes_{\mathcal O}E^0BQ
 \xrightarrow{\sim}E_{\bar k}^0BQ
\]

for symmetric groups and their products \cite[Section~3.23]{RezkCongruence}, compatibly with restrictions, transfers and Euler classes. Over $\bar k$, height classifies one-dimensional formal groups up to isomorphism \cite[Theorem~IV, pp.~270--271]{LazardFormalGroups}; choose an isomorphism with the height-$n$ Honda group. Applying the same coefficient comparison to the Honda theory over $\mathbb F_{p^n}$ identifies its extension with $E_{\bar k}$.

It follows that Strickland's identities hold after extension to $\mathcal O_{\bar k}$. Flat base change preserves the transfer images, restriction kernels, and annihilators, the latter being kernels of multiplication by finitely many generators. Faithful flatness descends these identities and the multiplication-by-$c_d^\Sigma$ isomorphism to $\mathcal O$.

Finally, if a change of oriented coordinate replaces $e_d$ by $e'_d=u_de_d$, with $u_d\in\mathcal O_d^\times$, then $\lambda'_d(b)=\lambda_d(u_db)$; thus perfectness is preserved. Formula~\eqref{eq:5.4} uses the Euler class and coordinate belonging to the same orientation.
\end{remark}

\begin{corollary}\label{cor:5.4}
There is a unique element $\delta_d\in \mathcal O_d$ such that
\begin{equation}\label{eq:5.9}
 \operatorname{Tr}_{\mathcal O_d/\mathcal O}(b)=\lambda_d(\delta_db)
 \qquad(b\in \mathcal O_d).
\end{equation}
It satisfies
\begin{equation}\label{eq:5.10}
 \boxed{e_d\delta_d=p^d.}
\end{equation}
Consequently the inverse different of $\mathcal O_d/\mathcal O$ is
\[
 \mathcal O_d^{\sharp}
 :=\{x\in \mathcal O_d[1/p]:\operatorname{Tr}_{\mathcal O_d/\mathcal O}(x\mathcal O_d)\subseteq\mathcal O\}
 =\frac{e_d}{p^d}\mathcal O_d.
\]
\end{corollary}

\begin{proof}
Existence and uniqueness follow from \eqref{eq:5.5}. Equations \eqref{eq:5.4} and \eqref{eq:5.9} give
\[
 \lambda_d((e_d\delta_d-p^d)b)
 =\operatorname{Tr}(e_db)-p^d\lambda_d(b)=0
 \qquad(b\in\mathcal O_d).
\]
Perfectness implies $e_d\delta_d=p^d$, so $\delta_d$ is invertible after inverting $p$. The $\lambda_d$-dual lattice is $\mathcal O_d$ itself; thus $\operatorname{Tr}(x\mathcal O_d)\subseteq\mathcal O$ is equivalent to $\delta_dx\in\mathcal O_d$. This gives $\mathcal O_d^\sharp=\delta_d^{-1}\mathcal O_d=(e_d/p^d)\mathcal O_d$.
\end{proof}

\subsection{The Euler norm}

Fix now a closed-point lift $z:\mathcal O\to W(k)$ and write
\[
 M_d=M_{d,z},\qquad e_{d,z}\in M_d
\]
for the image of $e_d$. The norm $\operatorname{Nm}_{M_d/W(k)}(b)$ means $\det(m_b)$. The relation \eqref{eq:5.10} shows that $e_{d,z}$ is invertible after inverting $p$, so its norm is nonzero.

\begin{lemma}\label{lem:5.5}
For $d\geq0$, with an empty sum when $d=0$,
\begin{equation}\label{eq:5.12}
 v_p\!\left(\operatorname{Nm}_{M_d/W(k)}(e_{d,z})\right)
 =\sum_{j=0}^{d-1}r_j.
\end{equation}
\end{lemma}

\begin{proof}
For $d=0$, use $e_0=1$. Suppose $d\geq1$ and put $G_z=z^*G$, as in Section~\ref{sec:3}.

The finite torsion scheme $G_z[p^d]$ has etale generic fibre. Choose a finite extension $K/W(k)[1/p]$ splitting it, with $v_p(p)=1$. The subgroup algebra then splits as
\[
 M_d\otimes_{W(k)}K
 \simeq
 \prod_{\substack{H\subset G_z[p^d](K)\\ |H|=p^d}}K.
\]
Multiplication is diagonal in this product, and determinants commute with scalar extension. Taking valuations gives
\[
 v_p(\operatorname{Nm}_{M_d/W(k)}(e_{d,z}))
 =
 \sum_{\substack{H\subset G_z[p^d](K)\\ |H|=p^d}}
 v_p(e_{d,z}(H)).
\]

For a fixed $H$, the corresponding transitive action is the regular action of $A_H=\mathbf Z_p^n/L_H$, with $A_H^\vee\simeq H$. Its reduced complex representation is
\[
 \mathbf C[A_H]-1
 \simeq
 \bigoplus_{0\ne h\in H}\mathbf C_h.
\]
In the chosen oriented coordinate $X$, the line indexed by $h$ has Euler value $X(h)$. Multiplicativity therefore gives
\begin{equation}\label{eq:5.13}
 e_{d,z}(H)
 =
 \prod_{\substack{h\in H\\ h\ne0}}X(h).
\end{equation}
This formula follows from the character decomposition and coordinate identification in \cite[Proposition~9.1 and its proof]{StricklandSymmetric}.

For $1\leq s\leq d$, the series $[p^s]_{G_z}(X)$ represents multiplication by $p^s$ on $G_z$. Its reduction modulo $p$ has order $p^{ns}$, so Weierstrass preparation gives
\[
 \frac{[p^s]_{G_z}(X)}{X}
 =
 U_s(X)W_s(X),
\]
where $U_s\in W(k)[[X]]^\times$ and $W_s\in W(k)[X]$ is monic distinguished of degree $p^{ns}-1$ \cite[Chapter~IV, \S~9, Theorem~9.2]{LangAlgebra}. Its roots are the coordinates of the nonzero $p^s$-torsion points. Since $W_s$ is distinguished, these coordinates have positive valuation.

The differential of $[p^s]$ at the identity is $p^s$. Thus $U_s(0)W_s(0)=p^s$, and $U_s(0)$ is a unit, giving
\begin{equation}\label{eq:5.14}
 v_p(W_s(0))=s.
\end{equation}

All torsion roots are simple: for translation $T_h$ by a $p^s$-torsion point, the identity $[p^s]\circ T_h=[p^s]$ gives
\[
 d[p^s]_h\circ d(T_h)_0=p^s\ne0.
\]
Since $T_h$ is an automorphism, $d(T_h)_0$ is invertible and $d[p^s]_h\ne0$. The nested root sets therefore give $W_{s-1}\mid W_s$ in $K[X]$, with quotient roots exactly the points of order $p^s$. Taking the valuation at $X=0$ and using \eqref{eq:5.14} yields
\begin{equation}\label{eq:5.15}
 \begin{aligned}
 \sum_{\substack{h\in G_z[p^s](K)\\
                 \operatorname{ord}(h)=p^s}}
 v_p(X(h))
 &=
 v_p(W_s(0))-v_p(W_{s-1}(0))\\
 &=s-(s-1)=1.
 \end{aligned}
\end{equation}
At $s=1$ use $W_0=1$. The calculation is independent of the oriented coordinate: if $X'=Xu(X)$ with $u(0)$ a unit, then positive valuation of $X(h)$ makes $u(X(h))$ a unit, so $v_p(X'(h))=v_p(X(h))$.

Identify the geometric $p$-power torsion with $\Lambda=(\mathbf Q_p/\mathbf Z_p)^n$. For a point $h$ of order $p^s$, quotient and inverse image give a bijection
\[
 \{H\subset\Lambda:|H|=p^d,\ h\in H\}
 \xrightarrow{\ \sim\ }
 \{\overline H\subset\Lambda/\langle h\rangle:
                  |\overline H|=p^{d-s}\}.
\]
Moreover,
\[
 (\Lambda/\langle h\rangle)^\vee
 =
 \{f\in\Lambda^\vee:f(h)=0\}
 \subseteq
 \Lambda^\vee\simeq\mathbf Z_p^n.
\]
This finite-index submodule is free of rank $n$, so duality gives $\Lambda/\langle h\rangle\simeq\Lambda$. Consequently $h$ belongs to exactly $r_{d-s}$ subgroups of order $p^d$.

Using \eqref{eq:5.13} and \eqref{eq:5.15}, reorder the finite sum to obtain
\[
 \begin{aligned}
 v_p\operatorname{Nm}_{M_d/W(k)}(e_{d,z})
 &=
 \sum_{\substack{H\subset\Lambda\\ |H|=p^d}}
 \ \sum_{0\ne h\in H}v_p(X(h))\\
 &=
 \sum_{s=1}^{d}
 r_{d-s}
 \sum_{\substack{h\in\Lambda\\
                  \operatorname{ord}(h)=p^s}}
 v_p(X(h))\\
 &=
 \sum_{s=1}^{d}r_{d-s}
 =\sum_{j=0}^{d-1}r_j.
 \end{aligned}
\]
\end{proof}

\subsection{Exact trace detection}

\begin{proposition}
\label{prop:5.6}
Let $W$ be a discrete valuation ring with uniformizer $p$, let $A$ be a $p$-torsion-free commutative $W$-algebra, and let $\tau:A\twoheadrightarrow B$ be a surjective $W$-algebra map with $B$ finite free over $W$. Fix $N\geq1$ and a descending filtration by ideals $A=F^0\supseteq\cdots\supseteq F^N$ satisfying
\begin{equation}\label{eq:5.16}
 F^{a+r}\cap p^rA=p^rF^a
 \qquad(a,r\geq0,\ a+r\leq N).
\end{equation}
For $0\leq a\leq N$ put
\[
 L^a(B)=\{b\in B:\operatorname{Tr}_{B/W}(bB)\subseteq p^aW\}.
\]
Suppose that
\[
 \tau(F^N)\subseteq L^N(B),\qquad
 \operatorname{length}_W(A/F^N)
 =\operatorname{length}_W(B/L^N(B))<\infty.
\]
Then
\[
 F^a=\tau^{-1}(L^a(B)),\qquad \tau(F^a)=L^a(B)
 \qquad(0\leq a\leq N).
\]
Moreover, $C=A/F^N$ has the perfect $W$-bilinear pairing
\begin{equation}\label{eq:5.19}
 C\times C\longrightarrow W/p^NW,\qquad
 (\bar x,\bar y)\longmapsto\operatorname{Tr}_{B/W}(\tau(x)\tau(y))\bmod p^N,
\end{equation}
and $\operatorname{Ann}_C(p^{N-a})=F^a/F^N$.
\end{proposition}

\begin{proof}
Each $L^a(B)$ is an ideal containing $p^aB$. The assumed inclusion induces a surjection $A/F^N\twoheadrightarrow B/L^N(B)$. Equality of the finite lengths makes it an isomorphism, so $F^N=\tau^{-1}(L^N(B))$.

For $0\leq a\leq N$ put $r=N-a$. Equation~\eqref{eq:5.16} and cancellation of $p^r$ in $A$ give $x\in F^a$ if and only if $p^rx\in F^N$. Thus
\[
 \begin{aligned}
 x\in F^a
 &\Longleftrightarrow p^rx\in F^N
 \Longleftrightarrow p^r\tau(x)\in L^N(B)\\
 &\Longleftrightarrow p^r\operatorname{Tr}_{B/W}(\tau(x)B)\subseteq p^NW
 \Longleftrightarrow \tau(x)\in L^a(B).
 \end{aligned}
\]
The final equivalence is cancellation in $W$. Surjectivity of $\tau$ gives the image identity. Passing the first equivalence to $C$ gives the asserted annihilator formula.

The equality $F^N=\tau^{-1}(L^N(B))$ makes \eqref{eq:5.19} well-defined and identifies its radical with zero. Equation~\eqref{eq:5.16} at $(a,r)=(0,N)$ gives $p^NC=0$. A finite-length module over $W$ is a direct sum of cyclic modules $W/p^eW$, here with $1\leq e\leq N$, and
\[
 \operatorname{Hom}_W(W/p^eW,W/p^NW)\simeq p^{N-e}W/p^NW
\]
has length $e$. Hence the injective adjoint $C\to\operatorname{Hom}_W(C,W/p^NW)$ has source and target of equal length and is an isomorphism. This proves perfectness without requiring $C$ to be free over $W/p^NW$.
\end{proof}

\begin{lemma}
\label{lem:5.7}
Let $B$ be a finite free commutative algebra over a discrete valuation ring $W$ with uniformizer $p$. Suppose that $\lambda:B\to W$ is a perfect Frobenius functional and that $e,\delta\in B$ satisfy
\[
 \operatorname{Tr}_{B/W}(b)=\lambda(\delta b)\quad(b\in B),
 \qquad e\delta=p^{N-1},\qquad N\geq1.
\]
Then
\[
 L^N(B)=peB,\qquad
 \operatorname{length}_W(B/L^N(B))
 =\operatorname{rank}_W(B)+v_p\operatorname{Nm}_{B/W}(e),
\]
where $v_p(p)=1$ and the norm is the determinant of multiplication by $e$.
\end{lemma}

\begin{proof}
Perfectness gives $b\in L^N(B)$ if and only if $\delta b\in p^NB$. If $\delta b=p^Nc$, multiplication by $e$ and cancellation of $p^{N-1}$ give $b=pec$. Conversely, $b=pec$ gives $\delta b=p^Nc$. This proves the lattice identity. Multiplication by $pe$ is injective, since $e\delta=p^{N-1}$ and $B$ is $p$-torsion free. Diagonalizing this endomorphism over the discrete valuation ring gives its cokernel as a sum of $W/p^{e_i}W$; it follows that its length is the valuation of its determinant. Finally, $\det(m_{pe})=p^{\operatorname{rank}_W(B)}\operatorname{Nm}_{B/W}(e)$.
\end{proof}

\begin{corollary}
\label{cor:5.8}
Under the hypotheses of Proposition~\ref{prop:5.6}, suppose also that $L^N(B)\subseteq pB$. Then
\[
 F^N=pF^{N-1}+\ker(\tau).
\]
\end{corollary}

\begin{proof}
Put $J=\ker(\tau)$. The equality $F^N=\tau^{-1}(L^N(B))$ gives $J\subseteq F^N$. If $y\in F^N$, write $\tau(y)=pb$ and lift $b$ to $x\in A$. Then $j=y-px\in J$, so
\[
 y-j\in F^N\cap pA=pF^{N-1}.
\]
This proves one inclusion; the reverse follows from \eqref{eq:5.16} and $J\subseteq F^N$. 
\end{proof}

Applying Rezk's result \cite[Theorem~7.1 and Corollary~7.2]{RezkWitt} to the perfect $k$-algebra $k$ and identifying $\mathbb W_E(k)$ with $\mathcal O$ via \eqref{eq:2.10} gives the following theorem. It supplies the intersection condition in Proposition~\ref{prop:5.6} with $N=d+1$.

\begin{theorem}\label{thm:5.9}
For every $a,r\geq0$,
\[
 I_{a+r}\cap p^r\mathcal O=p^rI_a.
\]
\end{theorem}

Multiplication by $p$ gives exact sequences
\[
 0\longrightarrow\mathcal O/I_j\xrightarrow{\,p\,}\mathcal O/I_{j+1}
 \longrightarrow(\mathcal O/p)/\overline I_{j+1}\longrightarrow0
 \qquad(j\geq0).
\]
Indeed, Theorem~\ref{thm:5.9} gives $pI_j\subseteq I_{j+1}$ and $I_{j+1}\cap p\mathcal O=pI_j$; cancellation of $p$ proves injectivity. The final term has $k$-dimension $r_j$ by Proposition~\ref{prop:3.1}, since its coefficient twist is an automorphism of $k$. Starting with $\mathcal O/I_0=0$, additivity of length gives
\begin{equation}\label{eq:5.23}
 \operatorname{length}_{W(k)}(\mathcal O/I_{d+1})=\sum_{j=0}^d r_j.
\end{equation}

For $0\leq a\leq d+1$ define the trace lattice
\[
 L_{a,d,z}
 =\{b\in M_d:\operatorname{Tr}_{M_d/W(k)}(bM_d)\subseteq p^aW(k)\}.
\]
Thus $b\in L_{a,d,z}$ means that its trace pairing with \emph{every} element of $M_d$ is divisible by $p^a$. If
\[
 M_d^\sharp=\{x\in M_d[1/p]:
\operatorname{Tr}(xM_d)\subseteq W(k)\},
\]
then the same condition can be written
\[
 L_{a,d,z}=M_d\cap p^aM_d^\sharp
 \qquad\text{inside }M_d[1/p].
\]
These are full $W(k)$-lattices: they lie between $p^aM_d$ and $M_d$. In particular $L_{0,d,z}=M_d$. At $d=0$, the algebra is $W(k)$ and the trace is the identity, so $L_{a,0,z}=p^aW(k)$ for $a=0,1$.

\begin{theorem}\label{thm:5.10}
For every $d\geq0$, every closed-point lift $z$, and every $0\leq a\leq d+1$,
\begin{equation}\label{eq:5.25}
 \boxed{\tau_{d,z}(I_a)=L_{a,d,z}.}
\end{equation}
This gives the characterization
\begin{equation}\label{eq:5.26}
 \boxed{
 I_a=
 \{x\in\mathcal O:
   \Theta_{d,z}(xy)\in p^aW(k)\text{ for every }y\in\mathcal O\}.
 }
\end{equation}
At index $d+1$ we have a more explicit identity
\begin{equation}\label{eq:5.27}
 \boxed{L_{d+1,d,z}=p e_{d,z}M_d.}
\end{equation}
Hence
\begin{equation}\label{eq:5.28}
 \boxed{
 \mathcal O/I_{d+1}\simeq M_d/(p e_{d,z}).
 }
\end{equation}
\end{theorem}

The isomorphism \eqref{eq:5.28} is induced by $\tau_{d,z}$.

\begin{proof}
Write $W=W(k)$. Base change of Proposition~\ref{prop:5.2} and Corollary~\ref{cor:5.4} gives a perfect functional $\lambda_{d,z}$ and elements $e_{d,z},\delta_{d,z}$ with $e_{d,z}\delta_{d,z}=p^d$. Lemma~\ref{lem:5.7}, with $N=d+1$, proves \eqref{eq:5.27}; together with Lemma~\ref{lem:5.5} it gives
\[
 \operatorname{length}_W(M_d/L_{d+1,d,z})
 =r_d+v_p\operatorname{Nm}_{M_d/W}(e_{d,z})
 =\sum_{j=0}^d r_j.
\]
This equals \eqref{eq:5.23}.

To prove the inclusion $\tau_{d,z}(I_{d+1})\subseteq L_{d+1,d,z}$, let $x\in I_{d+1}$ and $b\in M_d$. Proposition~\ref{prop:3.3} supplies $y\in\mathcal O$ with $\tau_{d,z}(y)=b$. Since $xy\in I_{d+1}$, the boundary case $a=d+1$, $r=d$ of Lemma~\ref{lem:4.4} gives
\[
 \operatorname{Tr}_{M_d/W}(\tau_{d,z}(x)b)=\Theta_{d,z}(xy)\in p^{d+1}W.
\]
Thus $\tau_{d,z}(I_{d+1})\subseteq L_{d+1,d,z}$.

To apply Proposition~\ref{prop:5.6}, take $A=\mathcal O$, $F^a=I_a$ and $N=d+1$. Temporarily give the ring $B=M_d$ the $W$-action $c\mathbin{\star}b=\sigma_d(c)b$, using source multiplication on the right. Then $\tau_{d,z}:A\to B$ is $W$-linear, and
\[
 \operatorname{Tr}_{B/W}(b)=\sigma_d^{-1}\!\left(\operatorname{Tr}_{M_d/W}(b)\right).
\]
Notice that, in a source basis the multiplication matrix for the transported structure is obtained by applying $\sigma_d^{-1}$ to each entry. This transport preserves finite freeness, rank, lengths and all trace lattices, since $\sigma_d$ is an automorphism fixing $p$. Theorem~\ref{thm:5.9} verifies \eqref{eq:5.16}, and the preceding paragraphs establish trace inclusion and equality of quotient lengths at index $d+1$. Proposition~\ref{prop:5.6} now gives \eqref{eq:5.25} and \eqref{eq:5.26}. Its case $a=d+1$ and \eqref{eq:5.27} give \eqref{eq:5.28}.
\end{proof}

Proposition~\ref{prop:5.6} also gives a pairing and an annihilator identity. Put $W=W(k)$ and $C_d=\mathcal O/I_{d+1}$. Its perfect $W$-bilinear pairing is
\[
 C_d\times C_d\longrightarrow W/p^{d+1}W,
 \qquad (x,y)\longmapsto
   \sigma_d^{-1}\!\left(\Theta_{d,z}(xy)\right)\bmod p^{d+1}.
\]
Applying $\sigma_d^{-1}$ compensates for the target coefficient action, as in Remark~\ref{rem:2.4}. For $0\leq a\leq d+1$, put $r=d+1-a$. The annihilator identity is
\[
 \operatorname{Ann}_{C_d}(p^r)=I_a/I_{d+1}.
\]

Write 
\[
M_d^*=\operatorname{Hom}_{W(k)}(M_d,W(k)).
\]
By an elementary divisor $p^j$ mean a diagonal entry $p^j$, up to a unit, in Smith normal form over the discrete valuation ring $W(k)$.

\begin{corollary}\label{cor:5.11}
The trace map
\[
 \mathcal T_{d,z}:M_d\longrightarrow M_d^*,
 \qquad
 b\longmapsto(c\mapsto\operatorname{Tr}_{M_d/W(k)}(bc))
\]
has elementary divisor $p^j$ with multiplicity
\[
 r_j-r_{j-1}
 \qquad(0\leq j\leq d).
\]
\end{corollary}

\begin{proof}
By Proposition~\ref{prop:5.2} and Corollary~\ref{cor:5.4}, the specialized trace pairing is the perfect $\lambda_{d,z}$-pairing composed with multiplication by $\delta_{d,z}$, where $e_{d,z}\delta_{d,z}=p^d$.
Since $M_d$ is $p$-torsion free, multiplication by $\delta_{d,z}$ is injective. Hence $\mathcal T_{d,z}$ is injective. Let
\[
 b_1\leq\cdots\leq b_{r_d}
\]
be its $p$-adic Smith exponents.

For $1\leq a\leq d+1$, Smith normal form gives
\[
 \operatorname{length}_{W(k)}(M_d/L_{a,d,z})
 =\sum_i\max(a-b_i,0).
\]
By Theorem~\ref{thm:5.10}, $\tau_{d,z}$ induces a $\sigma_d$-semilinear ring isomorphism $  \mathcal O/I_a\xrightarrow{\sim}M_d/L_{a,d,z}$, which preserves $W(k)$-length. Thus \eqref{eq:5.23}, with $d$ replaced by $a-1$, gives
\[
 \sum_i\max(a-b_i,0)
 =\sum_{j=0}^{a-1}r_j.
\]
For $a=0$ both sides are zero. Subtracting the identity for $a-1$ from that for $a$ yields
\[
 \#\{i:b_i\leq a-1\}=r_{a-1}.
\]
Taking $a=d+1$ shows that all $r_d$ exponents satisfy $b_i\leq d$. Therefore, for $0\leq j\leq d$,
\[
 \#\{i:b_i=j\}
 =r_j-r_{j-1},
\]
with $r_{-1}=0$.
\end{proof}

\section{Geometric recursion from trace detection}\label{sec:6}

Throughout this section, all $\mathcal O$-module constructions on $M_{d,z}$ use the target action $a\cdot m=\tau_{d,z}(a)m$, including Fitting ideals, tensor products and $pI_dM_{d,z}$. Traces still use the source $W(k)$-structure.

\begin{theorem}\label{thm:6.1}
For every finite height $n\geq1$, every $d\geq0$, and every lift $z$ of the closed point,
\[
 I_{d+1}=pI_d+J_{d,z}
 =pI_d+\operatorname{Fitt}^{\mathcal O}_0(M_{d,z}).
\]
\end{theorem}
\begin{proof}
Use $A=\mathcal O$, $F^a=I_a$, $N=d+1$, and $B=M_{d,z}$ with $W(k)$-action $c\star b=\sigma_d(c)b$. The proof of Theorem~\ref{thm:5.10} verifies the hypotheses of Proposition~\ref{prop:5.6}, and \eqref{eq:5.27} gives $L^N(B)=pe_{d,z}B\subseteq pB$. Corollary~\ref{cor:5.8} yields the first equality; Proposition~\ref{prop:3.3} gives the Fitting equality.
\end{proof}

\begin{restatedcorollaryB}
For every $d\geq0$ and every lift $z$ as above,
\[
 \boxed{I_{d+1}=pI_d+J_{d,z}.}
\]
Choose independently a closed-point lift $z_d:\mathcal O\to W(k)$ for every $d\geq0$. Define
\[
 G_0=\mathcal O,\qquad
 G_{d+1}=pG_d+\ker\!\left(
 \mathcal O\xrightarrow{t_d}\mathcal O_d\otimes_{s_d,\mathcal O,z_d}W(k)\right).
\]
Then $G_d=I_d$ for every $d$.
\end{restatedcorollaryB}

\begin{proof}[Proof of Corollary~\ref{cor:B}]
Starting with $G_0=I_0=\mathcal O$, apply Theorem~\ref{thm:6.1} successively with the independently chosen lifts $z_d$. If $G_d=I_d$, it gives $G_{d+1}=pI_d+J_{d,z_d}=I_{d+1}$.
\end{proof}

Write
\[
 \mathfrak T_{d,z}:=\operatorname{Spf}(\mathcal O/J_{d,z})\simeq\operatorname{Spf}(M_{d,z})
\]
for the target locus of the fixed source $G_z$, embedded in $\mathfrak X$ by $\tau_{d,z}$.

\begin{corollary}

For every $d\geq0$ and every closed-point lift $z$, there are canonical ring isomorphisms
\[
 \mathcal O/I_{d+1}
 \simeq M_{d,z}\otimes_{\mathcal O}\mathcal O/(pI_d)
 \simeq M_{d,z}/pI_dM_{d,z}.
\]
Inclusion and reduction modulo $p$ also induce canonical $\overline{\mathcal O}$-module isomorphisms
\[
 J_{d,z}/pJ_{d,z}
 \xrightarrow{\sim}I_{d+1}/pI_d
 \xrightarrow{\sim}\overline I_{d+1}.
\]
The maps are $j+pJ_{d,z}\mapsto j+pI_d$ and $x+pI_d\mapsto x+p\mathcal O$. Thus $pI_d+J_{d,z}$ is independent of $z$, although $J_{d,z}$ can vary.
\end{corollary}

\begin{proof}
Using $M_{d,z}\simeq\mathcal O/J_{d,z}$ and Theorem~\ref{thm:6.1}, the tensor product is $\mathcal O/(J_{d,z}+pI_d)=\mathcal O/I_{d+1}$. This quotient has finite $W(k)$-length by \eqref{eq:5.23}, so its quotient adic topology is discrete. Completion leaves it unchanged.

Since $J_{d,z}\subseteq I_{d+1}$ by Theorem~\ref{thm:6.1} and $I_{d+1}\subseteq I_d$ by the definition of the filtration in Section~\ref{sec:2}, inclusion defines the first module map. Theorem~\ref{thm:5.9} gives $I_{d+1}\cap p\mathcal O=pI_d$, therefore reduction defines the second isomorphism. Their composite is \eqref{eq:3.7}; hence the first map is an isomorphism.
\end{proof}

\section{The universal trace cokernel}\label{sec:7}

For Theorem~\ref{thm:C}, we lift compatible Smith forms from two discrete valuation rings to $\mathcal O$. After a faithfully flat extension at $(p)$, the subgroup algebra decomposes into torsors; their group orders determine trace exponents matching the specialized exponents of Corollary~\ref{cor:5.11}. This agreement permits the lift by Lemma~\ref{lem:7.4}. All module structures and base changes use the source map $s_d$.

Write $r_j^{(h)}$ for the number of order-$p^j$ subgroups of $(\mathbf Q_p/\mathbf Z_p)^h$, with $r_{-1}^{(h)}=0$. Thus $r_j=r_j^{(n)}$; at height zero use $r_0^{(0)}=1$ and $r_j^{(0)}=0$ for $j>0$.

Let $\mathcal G$ be the algebraic $p$-divisible group formed by the torsion schemes of $G$. At level $v\geq1$ its algebra is $\mathcal O[X]/(F_v(X))$, where $F_v$ is the distinguished Weierstrass polynomial of $[p^v](X)$. The algebra $\mathcal O_d$ represents rank-$p^d$ subgroups of $\mathcal G[p^d]$ after arbitrary algebraic base change \cite[Section~3.1, after Theorem~3.3]{StapletonSubgroups}. We base change these finite algebras to $\mathcal O_{(p)}$, retaining all their components. The reduction there is \emph{ordinary}: its connected part has height one.

For the multiplicative group $\mathbf G_m$, write $\mu_{p^v}=\ker([p^v]:\mathbf G_m\to\mathbf G_m)$ for $v\geq0$ and $\mu_{p^\infty}=\varinjlim_v\mu_{p^v}$. A torsor under a group scheme $D$ is a sheaf with a $D$-action which becomes $D$ acting by translation after a faithfully flat, locally finitely presented (fppf) cover. Underlined $\operatorname{Hom}$ denotes the fppf sheaf of group homomorphisms.

\begin{proposition}
\label{prop:7.1}
Suppose $n>1$. Put $U=\widehat{\mathcal O_{(p)}}$, and let $R$ be the completion of a strict henselization of $U$. After compatible trivializations of the etale factors, there is an exact sequence
\begin{equation}\label{eq:7.1}
 0\longrightarrow\mu_{p^\infty}\longrightarrow\mathcal G_R
 \xrightarrow{\ \pi\ }Q:=(\mathbf Q_p/\mathbf Z_p)^{n-1}\longrightarrow0.
\end{equation}
Fix $d\geq0$. Regard finite subgroups $K\subset Q$ as constant group schemes. For $|K|=p^j$, $0\leq j\leq d$, put $a=d-j$ and
\begin{equation}\label{eq:7.2}
 q_a:\mathcal G_R\longrightarrow\mathcal G':=\mathcal G_R/\mu_{p^a},
 \qquad \rho:\mathcal G'\longrightarrow Q,\qquad \rho q_a=\pi.
\end{equation}
The sheaf of group-homomorphism sections over $K$, defined for $R$-algebras $S$ by
\[
 X_K(S)=\{\sigma:K_S\longrightarrow\mathcal G'_S:
                 \rho_S\sigma=i_{K,S}\},
 \qquad i_{K,S}:K_S\hookrightarrow Q_S,
\]
is a torsor under the Cartier dual $K^D=\underline{\operatorname{Hom}}(K,\mathbf G_m)$.
It is represented by a finite locally free affine $R$-scheme of degree $p^j$.
Writing $C_K=\mathcal O(X_K)$, one has
\begin{equation}\label{eq:7.4}
 \mathcal O_d\otimes_{\mathcal O}R
 \simeq\prod_{j=0}^d\prod_{\substack{K\subset Q\\|K|=p^j}}C_K.
\end{equation}
No splitting of \eqref{eq:7.1} is assumed.
\end{proposition}

\begin{proof}
Choose auxiliary deformation parameters $v_1,\ldots,v_{n-1}$ and a coordinate $X$ such that $\mathcal O\simeq W(k)[[v_1,\ldots,v_{n-1}]]$ and $[p]_G(X)\equiv v_1X^p\pmod{(p,X^{p+1})}$, as in \cite[Section~2.1]{StapletonTranschromatic}. Since $n>1$, $v_1\notin(p)$, so the map to $U$ factors through
\[
 \mathcal O\longrightarrow L_1:=\widehat{\mathcal O[v_1^{-1}]}_{(p)}
 \longrightarrow U.
\]

Over $L_1$, the formal part has height one and the etale quotient has height $n-1$ \cite[Propositions~2.3, 2.4 and~2.7]{StapletonTranschromatic}; the same holds over $R$. The identity component is taken over the $p$-adic formal base. Its geometric special fibre is multiplicative by \cite[Theorem~IV]{LazardFormalGroups}, so the Cartier dual of each finite torsion level has etale special fibre. If $B$ is the finite locally free $R$-algebra of such a dual, base change of relative differentials gives $\Omega_{B/R}\otimes_R R/pR=0$. The module $\Omega_{B/R}$ is finite over $R$, so Nakayama's lemma gives $\Omega_{B/R}=0$. Finite flatness and finite presentation therefore make $B$ etale over $R$. The ring $R$ is a complete strictly henselian discrete valuation ring, and $\mathcal O_{(p)}\to R$ is faithfully flat with uniformizer $p$. Its finite etale group schemes are constant. Compatible trivializations of the quotient and the height-one dual give \eqref{eq:7.1}.

We first work over any $R/p^vR$-algebra $S$, with $v\geq1$. A rank-$p^d$ subgroup $H\subseteq\mathcal G_S[p^d]$ has open and closed fibres $H_b$ over $b\in Q[p^d]$, each finite locally free over $S$. Partition $\operatorname{Spec}(S)$ into open and closed pieces where their ranks are fixed. On each piece the positive-rank labels form a subgroup $K$: the corresponding fibres have sections fppf locally, and multiplication and inversion give closure. Translation by such a section identifies each fibre with $H_0$. Hence $H\to K_S$ is faithfully flat, $|K|=p^j$ and $\operatorname{rank}(H_0)=p^{d-j}$.

Put $a=d-j$. We claim that $H_0=\mu_{p^a,S}$, including for nonreduced $S$. Its Cartier dual $H_0^D$ is finite locally free. On a geometric fibre over an algebraically closed field $\kappa$ of characteristic $p$, the coordinate algebra of $\mu_{p^d}$ is $\kappa[T]/(T^{p^d})$, whose unique quotient of dimension $p^a$ is $\kappa[T]/(T^{p^a})$. Thus the subgroup on that fibre is $\mu_{p^a}$, and its dual is the etale cyclic group $\mathbf Z/p^a$. The relative differentials of $H_0^D$ therefore vanish on every geometric fibre, thus vanish by local Nakayama; finite flatness makes $H_0^D$ etale over $S$. The dual inclusion gives a map $(\mathbf Z/p^d)_S\to H_0^D$ that is surjective on geometric fibres. Etale locally both groups are constant, so this map is finite etale and surjective. Its kernel is an open and closed subgroup of $(\mathbf Z/p^d)_S$ of rank $p^{d-a}$. The unique subgroup of that order in $\mathbf Z/p^d$ is $p^a(\mathbf Z/p^d)$, so the kernel is $(p^a(\mathbf Z/p^d))_S$. Dualizing this standard cyclic quotient proves the claimed equality of embedded subgroup schemes.

Under \eqref{eq:7.2}, $H/\mu_{p^a,S}$ gives a section over $K_S$. Conversely, a section $\sigma$ gives $H=q_a^{-1}(\sigma(K_S))$, finite locally free of rank $p^a|K|=p^d$. Its $p^j$-multiple lies in $\mu_{p^a,S}$, so it is killed by $p^d$. These constructions are inverse and commute with base change.

For the torsor assertion, let $S$ be an arbitrary $R$-algebra. Multiplication by $p^a$ identifies $\ker(\rho)=\mu_{p^\infty}/\mu_{p^a}$ with $\mu_{p^\infty}$. Two sections over $K_S$ differ by an element of
\[
 \underline{\operatorname{Hom}}(K,\mu_{p^\infty})
 =\underline{\operatorname{Hom}}(K,\mathbf G_m)=K^D.
\]
Sections exist fppf locally. Write $K=\bigoplus_i\mathbf Z/p^{\lambda_i}$. Since $\rho$ is surjective as an fppf sheaf, after an fppf cover of $\operatorname{Spec}(S)$ each generator lifts to a section $g_i$ of $\mathcal G'$. Its multiple $p^{\lambda_i}g_i$ lies in the kernel. Multiplication by $p^{\lambda_i}$ on $\mu_{p^\infty}$ is fppf surjective, therefore after a further cover choose a kernel section $h_i$ with the same multiple. There are finitely many generators, so common fppf cover suffices for these choices. The lifts $g_i-h_i$ satisfy the cyclic relations and, by commutativity, define a group-homomorphism section over $K$. It follows that $X_K$ is a $K^D$-torsor over $R$. Effective affine descent and descent of finite local freeness represent it by an affine scheme of degree $|K|$ \cite[Lemmas~35.37.1 and~35.23.32]{StacksDescent}.

The correspondence over all $R/p^vR$-algebras, including nonreduced ones, gives compatible isomorphisms of representing algebras modulo $p^v$. Both sides of \eqref{eq:7.4} are finite $R$-modules, thus $p$-adically complete; taking inverse limits proves the formula. For $K=0$, the section is unique and $C_K=R$.
\end{proof}

\begin{lemma}\label{lem:7.2}
Let $R$ be a commutative local ring, $K$ a finite abelian group of order $m$, and $C$ the algebra of an fppf torsor under $K^D=\operatorname{Spec} R[K]$, where $R[K]$ is the group algebra. There is a grading
\[
 C=\bigoplus_{\chi\in K}C_\chi,\qquad C_0=R,\qquad
 C_\chi\otimes_R C_\psi\xrightarrow{\sim}C_{\chi+\psi},
\]
by invertible (rank-one locally free) $R$-modules, with the last map given by multiplication. For the degree-zero projection $\pi_0:C\to R$,
\begin{equation}\label{eq:7.8}
 \beta(x,y):=\pi_0(xy)\ \text{is perfect},\qquad
 \operatorname{Tr}_{C/R}(xy)=m\beta(x,y).
\end{equation}
When $R$ is a discrete valuation ring with uniformizer $p$ and $m=p^j$, all $m$ Smith exponents of its trace map are $j$.
\end{lemma}

\begin{proof}
Write the coaction as $\alpha(c)=\sum_\chi c_\chi\otimes z^\chi$, where $z^\chi$ is the group-algebra basis, with comultiplication $z^\chi\mapsto z^\chi\otimes z^\chi$. Counit and coassociativity give $c=\sum_\chi c_\chi$ and $\alpha(c_\chi)=c_\chi\otimes z^\chi$. Coefficient uniqueness gives the direct-sum grading, and multiplicativity of $\alpha$ makes it an algebra grading.

After a faithfully flat extension trivializing the torsor, this is the usual grading of $R[K]$. Descent makes $C$ finite locally free of rank $m$. Its graded summands are projective, and their formation commutes with base change. Their rank-one property and the isomorphisms $R\to C_0$ and $C_\chi\otimes_R C_\psi\to C_{\chi+\psi}$ therefore descend.

Multiplication by $x\in C_\chi$ sends $C_\psi$ to $C_{\chi+\psi}$. For $\chi\ne0$ its matrix in homogeneous bases has zero diagonal; for $x\in C_0=R$ its trace is $mx$. Thus $\operatorname{Tr}_{C/R}=m\pi_0$. The isomorphisms $C_\chi\otimes_R C_{-\chi}\to R$ make $\beta$ perfect, proving \eqref{eq:7.8} without division by $m$. Over a DVR the trace matrix is $p^j$ times an invertible matrix.
\end{proof}

\begin{corollary}
\label{cor:7.3}
Over $\mathcal O_{(p)}$, the trace map of $\mathcal O_d$ has Smith exponent $j$ with multiplicity
\begin{equation}\label{eq:7.9}
 p^j r_j^{(n-1)}=r_j^{(n)}-r_{j-1}^{(n)}
 \qquad(0\leq j\leq d).
\end{equation}
\end{corollary}

\begin{proof}
  For $n=1$ one has $\mathcal O_d=\mathcal O$ and trace matrix $(1)$, which gives the formula. For $n>1$, Proposition~\ref{prop:7.1} and Lemma~\ref{lem:7.2} give exponent $j$ repeated $p^j$ times for each of the $r_j^{(n-1)}$ torsors. Put $F_h(T)=\sum_{j\geq0}r_j^{(h)}T^j$. The rank series \eqref{eq:5.1} gives $(1-T)F_n(T)=F_{n-1}(pT)$; comparison of coefficients proves \eqref{eq:7.9}. The faithfully flat extension $\mathcal O_{(p)}\to R$ preserves the uniformizer $p$ and thus the valuations of all determinantal ideals. Their successive differences recover the same Smith list over $\mathcal O_{(p)}$.
\end{proof}

The specialized Smith exponents of Corollary~\ref{cor:5.11} do not give a diagonal form over $\mathcal O$. Indeed, consider the matrix over $A=W(k)[[u]]$
\[
 T=\begin{pmatrix}p&u\\0&p\end{pmatrix}.
\]
Its specialization at $u=0$ is $\operatorname{diag}(p,p)$, with Smith exponents $(1,1)$. Over $A_{(p)}$, however, $u$ is a unit, so the exponents are $(0,2)$. Thus $T$ and $\operatorname{diag}(p,p)$ have the same determinant $p^2$ and the same specialized Smith form, but are not equivalent over $A$: such an equivalence would persist after localization. For the universal trace map, Corollary~\ref{cor:7.3} supplies the localized Smith list, which agrees with the specialized list of Corollary~\ref{cor:5.11}. The following lemma lifts this common list to a diagonal form over $\mathcal O$.

\begin{lemma}\label{lem:7.4}
Put
\[
 A=W(k)[[u_1,\ldots,u_{n-1}]].
\]
Let $T$ be an $N\times N$ matrix over $A$, with $N\geq0$. Suppose its Smith exponents over $A_{(p)}$ and those of its specialization $u_1=\cdots=u_{n-1}=0$ over $W(k)$ are the same list of $N$ finite nonnegative integers $e_1\leq\cdots\leq e_N$. Then there exist $P,Q\in\operatorname{GL}_N(A)$ such that
\[
 PTQ=\operatorname{diag}(p^{e_1},\ldots,p^{e_N}).
\]
\end{lemma}

\begin{proof}
The case $N=0$ is empty. For $N>0$, first observe that
\begin{equation}\label{eq:7.11}
 p^eA_{(p)}\cap A=p^eA\qquad(e\geq0).
\end{equation}
For $e=0$ this is immediate. Otherwise write $sf=p^eg$ with $s\notin(p)$. Since $A/pA=k[[u_1,\ldots,u_{n-1}]]$ is a domain, reduction modulo $p$ gives $f=pf_1$. Cancel $p$ and repeat to obtain $f\in p^eA$.

The smallest localized exponent and \eqref{eq:7.11} show that $U=p^{-e_1}T$ has entries in $A$. Since the specialized smallest exponent is also $e_1$, some entry of $U$ specializes to a unit of $W(k)$, hence is a unit of $A$.

We move this entry to the upper-left corner and normalize its unit factor. Every entry of $T$ is divisible by the resulting pivot $p^{e_1}$, so integral row and column operations give
\[
 \begin{pmatrix}p^{e_1}&0\\0&T'\end{pmatrix}.
\]
The operations remain invertible under both base changes. All entries of $T'$ remain divisible by $p^{e_1}$, so its Smith list at each DVR is the original list with one occurrence of $e_1$ removed. Induction proves the claim.
\end{proof}

The \emph{discriminant ideal} of a finite free algebra is generated by the determinant of its trace Gram matrix $(\operatorname{Tr}(b_i b_j))$. A change of basis multiplies this determinant by a unit, so the ideal is intrinsic.

\begin{theorem}\label{thm:7.5}
For $n\geq1$ and $d\geq0$, put $\mathcal O_d^*=\operatorname{Hom}_{\mathcal O,s_d}(\mathcal O_d,\mathcal O)$. The universal source trace map
\[
 \mathcal T_d:\mathcal O_d\longrightarrow \mathcal O_d^*,
 \qquad
 b\longmapsto\bigl(c\mapsto\operatorname{Tr}_{(\mathcal O_d,s_d)/\mathcal O}(bc)\bigr),
\]
admits independent changes of domain and codomain bases, given by $P,Q\in\operatorname{GL}_{r_d}(\mathcal O)$, such that
\begin{equation}\label{eq:7.12}
 P[\mathcal T_d]Q=
 \boxed{
 \bigoplus_{j=0}^{d}
 p^j\operatorname{id}_{\mathcal O^{\,r_j-r_{j-1}}}
 }.
\end{equation}
We therefore obtain, noncanonically as an $\mathcal O$-module,
\begin{equation}\label{eq:7.13}
 \boxed{
 \operatorname{coker}(\mathcal T_d)
 \simeq
 \bigoplus_{j=1}^{d}
 (\mathcal O/p^j\mathcal O)^{\,r_j-r_{j-1}}.
 }
\end{equation}
The discriminant ideal is
\begin{equation}\label{eq:7.14}
 \boxed{
 \operatorname{disc}(\mathcal O_d/\mathcal O)
 =\left(p^{\,d r_d-\sum_{j=0}^{d-1}r_j}\right).
 }
\end{equation}
\end{theorem}

\begin{proof}
Finite freeness makes the source dual and trace pairing commute with base change: apply the coefficient map to the multiplication matrices. Corollary~\ref{cor:7.3} gives the Smith list over $\mathcal O_{(p)}$, and Corollary~\ref{cor:5.11} gives the same list at the source lift $u_i=0$. Lemma~\ref{lem:7.4} proves \eqref{eq:7.12}. Taking its cokernel gives \eqref{eq:7.13}, since the exponent-zero block is invertible. Summing the exponents gives
\[
 \sum_{j=1}^{d}j(r_j-r_{j-1})
 =d r_d-\sum_{j=0}^{d-1}r_j,
\]
which proves \eqref{eq:7.14}. At $d=0$ the cokernel is zero and the discriminant is the unit ideal.
\end{proof}

\begin{remark}
The decomposition \eqref{eq:7.13} of source $\mathcal O$-modules is noncanonical. The trace map and its Fitting ideals are canonical: for $0\leq i<r_d$, $\operatorname{Fitt}_i^\mathcal O(\operatorname{coker}\mathcal T_d)$ is generated by the $(r_d-i)$-minors of $[\mathcal T_d]$, and for $i\geq r_d$ it equals $\mathcal O$.
\end{remark}

\begin{remark}[Example: height two, $d=1$]\label{rem:7.7}
At height two, $Q=\mathbf Q_p/\mathbf Z_p$ has one subgroup of each order $p^j$. For $d=1$, \eqref{eq:7.4} becomes
\[
 \mathcal O_1\otimes_\mathcal O R\simeq R\times C,
 \qquad \operatorname{Spec}(C)\text{ a }\mu_p\text{-torsor over }R.
\]
The trace blocks are $(1)$ and $p$ times an invertible $p\times p$ matrix by Lemma~\ref{lem:7.2}. With $r_0=1$, $r_1=p+1$, Theorem~\ref{thm:7.5} gives
\[
 \operatorname{coker}(\mathcal T_1)\simeq(\mathcal O/p\mathcal O)^p,
 \qquad
 \operatorname{disc}(\mathcal O_1/\mathcal O)=(p^p).
\]

For $p=2$, use the model in Section~\ref{sec:3.1} over $\mathcal O=\mathbf Z_2[[u]]$. In the source basis $1,b,b^2$, multiplication by $b$ has matrix
\[
 D=\begin{pmatrix}0&0&2\\1&0&u\\0&1&0\end{pmatrix}.
\]
The relation $b^3=ub+2$ gives
$\operatorname{Tr}(1)=3$, $\operatorname{Tr}(b)=0$, $\operatorname{Tr}(b^2)=2u$, $\operatorname{Tr}(b^3)=6$ and
$\operatorname{Tr}(b^4)=2u^2$. Thus the trace Gram matrix is
\begin{equation}\label{eq:7.15}
 G=\begin{pmatrix}
 3&0&2u\\
 0&2u&6\\
 2u&6&2u^2
 \end{pmatrix},
 \qquad \det G=4(u^3-27).
\end{equation}
Splitting the unit pivot $3$ by row and column operations leaves
\[
 2\begin{pmatrix}u&3\\3&u^2/3\end{pmatrix}.
\]
The inner determinant $(u^3-27)/3$ is a unit of $\mathcal O$. Independent row and column operations therefore give $\operatorname{diag}(1,2,2)$ and cokernel $(\mathcal O/2)^2$ over the full deformation ring.
\end{remark}

\section{Trace operations and their composition}\label{sec:8}

Rezk's integral Hecke action sums operations associated with quotient isogenies of prescribed kernel type \cite[Section~14.2 and Proposition~14.3]{RezkLogarithm}; see \cite[Section~2.4 and Proposition~4.5]{ZhuHecke} for explicit height-two calculations. We study integral weighted traces on subgroup algebras and their composition.

For $d\geq0$, the source dual $\mathcal A[d]=\operatorname{Hom}_{\mathcal O,s_d}(\mathcal O_d,\mathcal O)$ represents additive operations of weight $p^d$ on degree-zero homotopy of $K(n)$-local commutative Morava $E$-algebras \cite[Section~3.12 and proof of Proposition~3.21]{RezkKoszul}, using the subgroup identification of Section~\ref{sec:2.6}. The weight records the subgroup order. Put $\mathcal A=\bigoplus_{d\geq0}\mathcal A[d]$, the algebraic direct sum graded by $d$. Its $\mathfrak m$-adic completion as a source module is the full ring of natural additive operations \cite[Section~3.22 and Proposition~3.23]{RezkKoszul}; our results concern the direct sum.

On $\mathcal O$ we write $x\gamma=\gamma(t_d(x))$. The product $\gamma\eta$ denotes ordinary composition of the natural operations represented by $\gamma$ and $\eta$, with $\eta$ applied first. Postfix notation reverses the displayed order of application: $x(\gamma\eta)=(x\eta)\gamma$.

Give $\mathcal A[d]$ its $\mathcal O_d$-module structure 
\[
 (b\lambda)(c)=\lambda(bc),
 \qquad b,c\in \mathcal O_d,\quad\lambda\in\mathcal A[d],
\]
and write the submodule of trace operations as
\[
 \mathcal A_{\mathrm{tr}}[d]
 =\{c\mapsto\operatorname{Tr}_{\mathcal O_d/\mathcal O}(bc):b\in \mathcal O_d\}\subseteq\mathcal A[d].
\]
This is an $\mathcal O_d$-submodule.

\begin{corollary}\label{cor:8.1}
There is a canonical exact sequence of source $\mathcal O$-modules
\[
 0\longrightarrow \mathcal O_d\xrightarrow{\mathcal T_d}\mathcal A[d]
 \longrightarrow C_d^{\mathrm{tr}}\longrightarrow0,
 \qquad C_d^{\mathrm{tr}}=\mathcal A[d]/\mathcal A_{\mathrm{tr}}[d],
\]
and a noncanonical decomposition
\begin{equation}\label{eq:8.1}
 C_d^{\mathrm{tr}}\simeq
 \bigoplus_{j=1}^{d}(\mathcal O/p^j)^{r_j-r_{j-1}}.
\end{equation}
Under the isomorphism $b\mapsto\lambda_d(b-)$ of Proposition~\ref{prop:5.2}, the submodule $\mathcal A_{\mathrm{tr}}[d]$ corresponds to $\delta_d\mathcal O_d$, so that $C_d^{\mathrm{tr}}\simeq \mathcal O_d/(\delta_d)$ as $\mathcal O_d$-modules. For $n\geq2$ and $d\geq1$, the least integer $s\geq0$ with $p^s\mathcal A[d]\subseteq\mathcal A_{\mathrm{tr}}[d]$ is $d$. At height one, and also for $d=0$, the cokernel vanishes.
\end{corollary}

\begin{proof}
For $b,c\in\mathcal O_d$, equation~\eqref{eq:5.9} gives $\mathcal T_d(b)(c)=\operatorname{Tr}(bc)=\lambda_d(\delta_dbc)$. Thus Proposition~\ref{prop:5.2} identifies $\mathcal T_d$ with multiplication by $\delta_d$. Since $e_d\delta_d=p^d$ and $\mathcal O_d$ is $p$-torsion free, this map is injective. Theorem~\ref{thm:7.5} gives \eqref{eq:8.1}. For $\gamma=\lambda_d(b-)$,
\[
 p^d\gamma=\operatorname{Tr}_{\mathcal O_d/\mathcal O}(e_db-).
\]
The exponent is sharp since \eqref{eq:7.9} gives $r_d-r_{d-1}>0$ for $n\geq2,d\geq1$. At height one all $r_j=1$, and at $d=0$ the sum is empty.
\end{proof}

Thus $\mathcal A_{\mathrm{tr}}[d]$ is a finite source $\mathcal O$-submodule containing $p^d\mathcal A[d]$ and becomes all of $\mathcal A[d]$ after inverting $p$.

\begin{remark}[The running example as trace operations]\label{rem:8.2}
In the model \eqref{eq:3.9}, the rows of the trace matrix \eqref{eq:7.15}, expressed in the basis dual to $1,b,b^2$, give $\operatorname{Tr}(-)$, $\operatorname{Tr}(b-)$ and $\operatorname{Tr}(b^2-)$. The first is Rezk's trace operator \cite[Section~2.8]{RezkHeightTwo}. These three functionals generate $\mathcal A_{\mathrm{tr}}[1]$, whose quotient in $\mathcal A[1]$ is $(\mathcal O/2)^2$ by Remark~\ref{rem:7.7}.
\end{remark}

The tensor products below use the source and target coefficient actions
\[
 (c\gamma)(b)=c\gamma(b),\qquad
 (\gamma c)(b)=\gamma(t_d(c)b)
 \quad(c\in\mathcal O,\ \gamma\in\mathcal A[d],\ b\in \mathcal O_d).
\]
The first multiplies the output by $c$; the second inserts multiplication by $c$ before the operation. Reversing products gives Rezk's right-operation convention \cite[Section~5.7]{RezkCongruence}.

A two-step isogeny flag is a quotient of degree $p^a$ followed by one of degree $p^b$. Matching the first target with the second source gives its coordinate algebra
\[
 C_{a,b}=\mathcal O_a\widehat\otimes_{t_a,\mathcal O,s_b}\mathcal O_b.
\]
It can also be described as a chain $H_a\subseteq H_{a+b}\subseteq G$ with $|H_a|=p^a$ and $|H_{a+b}/H_a|=p^b$. The tensor relation
\[
 x\,t_a(c)\otimes y=x\otimes s_b(c)y
 \qquad(c\in\mathcal O)
\]
shows the equivalence of the two descriptions.

Choose an $\mathcal O$-basis $y_1,\ldots,y_{r_b}$ of $\mathcal O_b$ for the source action $s_b$. The ordinary tensor product $C^{\mathrm{alg}}=\mathcal O_a\otimes_{t_a,\mathcal O,s_b}\mathcal O_b$ is free over $\mathcal O_a$ on $1\otimes y_j$, hence complete. If $\mathfrak n_a,\mathfrak n_b$ are the respective maximal ideals, the Artinian closed fibre gives $\mathfrak n_b^N\subseteq s_b(\mathfrak m)\mathcal O_b$ for some $N$. The image of $\mathfrak n_b^N$ in $C^{\mathrm{alg}}$ therefore lies in $\mathfrak n_a C^{\mathrm{alg}}$, since $t_a$ is local. Thus $\mathfrak n_a C^{\mathrm{alg}}$ and $\mathfrak n_a C^{\mathrm{alg}}+\mathfrak n_b C^{\mathrm{alg}}$ define equivalent adic topologies, so completion adds no elements:
\begin{equation}\label{eq:8.3}
 C_{a,b}=\mathcal O_a\otimes_{t_a,\mathcal O,s_b}\mathcal O_b,\qquad
 \operatorname{rank}_{\mathcal O,s_a}C_{a,b}=r_a r_b.
\end{equation}
In particular $C_{a,b}$ is $p$-torsion free, and every element is a finite sum of pure tensors.

The map that forgets $H_a$ in the chain, called the \emph{flag-forgetting map}, composes the two isogenies. Its pullback on coordinate rings is
\[
 \phi_{a,b}=\operatorname{comp}_{a,b}:
 \mathcal O_{a+b}\longrightarrow C_{a,b},
\]
with
\[
 \operatorname{comp}_{a,b}\circ s_{a+b}=s_a\otimes1,\qquad
 \operatorname{comp}_{a,b}\circ t_{a+b}=1\otimes t_b.
\]
Under the identification of Section~\ref{sec:2.6}, this map is dual to composition of power operations \cite[Proposition~3.41(3)]{BSY}. We spell out the duality to fix the multiplication order.

Use the right coefficient action through $t_a$ on $\mathcal A[a]$ and the left coefficient action through $s_b$ on $\mathcal A[b]$. Put
\[
 C_{a,b}^{\vee}=\operatorname{Hom}_{\mathcal O}(C_{a,b},\mathcal O),
\]
where $\mathcal O$ acts on $C_{a,b}$ through the source map $\mathcal O\xrightarrow{s_a}\mathcal O_a\to C_{a,b}$. For $\gamma\in\mathcal A[a]$ and $\eta\in\mathcal A[b]$, define
\[
 \langle\gamma,\eta\rangle(x\otimes y)
 =\gamma\bigl(x\,t_a(\eta(y))\bigr).
\]
Source linearity of $\eta$ gives
\[
 \begin{aligned}
 \langle\gamma,\eta\rangle(x\otimes s_b(c)y)
 &=\gamma\bigl(x\,t_a(c\eta(y))\bigr)\\
 &=\langle\gamma,\eta\rangle(x t_a(c)\otimes y).
 \end{aligned}
\]
This descends to a source-linear functional on $C_{a,b}$. The coefficient actions on the operation modules satisfy
\begin{equation}\label{eq:8.6}
 \langle\gamma\cdot c,\eta\rangle(x\otimes y)
 =\langle\gamma,c\eta\rangle(x\otimes y).
\end{equation}
Hence the pairing induces a map
\[
 \mathcal A[a]_{\mathrm t}\otimes_{\mathcal O}{}_{\mathrm s}\mathcal A[b]
 \longrightarrow C_{a,b}^{\vee}.
\]
Here the subscripts record these right and left coefficient actions.

This map is an isomorphism. To see its inverse, let $y_j^*$ be the source dual basis of $\mathcal O_b$. For $f\in C_{a,b}^{\vee}$, set
\[
 \gamma_j(x)=f(x\otimes y_j),\qquad
 f\longmapsto\sum_{j=1}^{r_b}\gamma_j\otimes y_j^*.
\]
Each $\gamma_j$ is source-linear. Writing $y=\sum_j s_b(y_j^*(y))y_j$ shows
\[
 f(x\otimes y)
 =\sum_j\gamma_j\bigl(x\,t_a(y_j^*(y))\bigr).
\]
Conversely, expanding $\eta=\sum_j\eta(y_j)y_j^*$ and using \eqref{eq:8.6} recovers $\gamma\otimes\eta$. The multiplication dual to composition of isogenies is
\begin{equation}\label{eq:8.8}
 \gamma\eta=\langle\gamma,\eta\rangle\circ\phi_{a,b}.
\end{equation}
Evaluation on $t_{a+b}(c)$ gives $\gamma(t_a(\eta(t_b(c))))$, confirming this order on coefficients; compare \cite[Sections~5.7, 11.12 and~12.4]{RezkCongruence}.

For the universal deformation $G/\mathcal O$, let $\mathfrak a_G$ be the ideal of functions vanishing at the identity. Its cotangent line
\[
 \omega_G=\mathfrak a_G/\mathfrak a_G^2
\]
is canonically the module of invariant differentials \cite[Section~11.17]{RezkCongruence}. In a coordinate $X$, the class of $X$ corresponds to $dX$ at the identity. This module is free of rank one; choose a generator $\omega$.

For the universal quotient $q_d:s_d^*G\to t_d^*G$ of Section~\ref{sec:2.6}, pullback is an $\mathcal O_d$-linear map $t_d^*\omega_G\to s_d^*\omega_G$. Its \emph{differential multiplier} $h_d\in\mathcal O_d$ is defined by
\begin{equation}\label{eq:8.9}
 q_d^*(t_d^*\omega)=h_d\,s_d^*\omega.
\end{equation}

\begin{lemma}
There are units $\xi_d\in \mathcal O_d^\times$ and elements $\kappa_d\in \mathcal O_d$ with
\[
 h_d=\xi_de_d,\qquad h_d\kappa_d=p^d,\qquad h_0=\kappa_0=1.
\]
The elements $h_d$ and $\kappa_d$ satisfy
\begin{equation}\label{eq:8.11}
 \phi_{a,b}(h_{a+b})=h_a\otimes h_b,\qquad
 \phi_{a,b}(\kappa_{a+b})=\kappa_a\otimes\kappa_b.
\end{equation}
\end{lemma}

\begin{proof}
Let $f_H(X)$ be the monic equation cutting out the universal subgroup. Its linear coefficient equals $e_d$ up to a unit \cite[proof of Theorem~9.2]{StricklandSymmetric}, with coefficient extension as in Remark~\ref{rem:5.3}. A quotient coordinate pulls back to $v(X)f_H(X)$ with $v(X)$ a unit \cite[Theorem~5.3 and its proof]{StricklandSubgroups}. Since $H$ contains the identity section, $f_H(0)=0$, so the derivative at $X=0$ is $v(0)f_H'(0)$. Changing the chosen generators of the modules of invariant differentials introduces only units, so $h_d=\xi_de_d$. Hence
\[
 \kappa_d:=\xi_d^{-1}\delta_d\in \mathcal O_d,\qquad h_d\kappa_d=p^d
\]
by Corollary~\ref{cor:5.4}.

A subgroup of order $p^d$ is killed by $p^d$ \cite[Proposition~5.2]{StricklandSubgroups}. The quotient property therefore gives a unique map $q_d^{\mathrm c}$ with
\[
 s_d^*G\xrightarrow{q_d}t_d^*G
 \xrightarrow{q_d^{\mathrm c}}s_d^*G,
 \qquad q_d^{\mathrm c}\circ q_d=[p^d].
\]
On the closed fibre, $[p^d]$ and $q_d$ have degrees $p^{nd}$ and $p^d$, so $q_d^{\mathrm c}$ has degree $p^{(n-1)d}$. It is therefore finite flat \cite[Proposition~3.3 and its proof]{StricklandSubgroups}; we call it the \emph{complementary isogeny} to $q_d$. Write $(q_d^{\mathrm c})^*(s_d^*\omega)=\widetilde\kappa_d\,t_d^*\omega$. Since $[p^d](X)=p^dX+O(X^2)$, pullback gives
\[
 h_d\widetilde\kappa_d\,s_d^*\omega
 =q_d^*(q_d^{\mathrm c})^*(s_d^*\omega)
 =[p^d]^*(s_d^*\omega)=p^d s_d^*\omega.
\]
Thus $h_d\widetilde\kappa_d=h_d\kappa_d=p^d$. Multiplying their difference by $\kappa_d$ gives $p^d(\widetilde\kappa_d-\kappa_d)=0$; $p$-torsion-freeness implies $\widetilde\kappa_d=\kappa_d$, hence
\[
 (q_d^{\mathrm c})^*(s_d^*\omega)=\kappa_d\,t_d^*\omega.
\]
Over $C=C_{a,b}$ write $i_a(x)=x\otimes1$ and $i_b(y)=1\otimes y$, and name the three coefficient maps
\[
 s=i_a s_a,\qquad \mu=i_a t_a=i_b s_b,\qquad t=i_b t_b.
\]
They give three formal groups $G_{\mathrm s}=s^*G$, $G_{\mathrm m}=\mu^*G$ and $G_{\mathrm t}=t^*G$, and the two base-changed quotient maps form the chain
\[
 \begin{tikzcd}[column sep=large]
 G_{\mathrm s}\arrow[r,"Q_a=i_a^*q_a"]
 &G_{\mathrm m}\arrow[r,"Q_b=i_b^*q_b"]
 &G_{\mathrm t}.
 \end{tikzcd}
\]
The composite is $\phi_{a,b}^*q_{a+b}$ under the unique marking-preserving target identification; uniqueness also gives coherence for further composition \cite[Sections~11.3-11.6]{RezkCongruence}. Thus
\[
 Q_b^*(t^*\omega)=(1\otimes h_b)\mu^*\omega,\qquad
 Q_a^*(\mu^*\omega)=(h_a\otimes1)s^*\omega.
\]
The chain rule proves the first identity in \eqref{eq:8.11}. Both $\phi_{a,b}(\kappa_{a+b})$ and $\kappa_a\otimes\kappa_b$ have product $p^{a+b}$ with $h_a\otimes h_b$. Since $(h_a\otimes h_b)(\kappa_a\otimes\kappa_b)=p^{a+b}$ and $C_{a,b}$ is $p$-torsion free by \eqref{eq:8.3}, multiplication by $h_a\otimes h_b$ is injective. Cancellation proves the second identity. At $d=0$ both isogenies are the identity, giving $h_0=\kappa_0=1$.
\end{proof}

The identities \eqref{eq:8.11} use the chosen invariant differential to fix the units in the comparison with the Euler classes. The following theorem shows that precomposition with multiplication by $\kappa_d$ has image $\mathcal A_{\mathrm{tr}}[d]$ and respects composition of operations.

\begin{theorem}
\label{thm:8.4}
For $\gamma\in\mathcal A[d]$ define
\begin{equation}\label{eq:8.15}
 (\mathcal V_d^\omega\gamma)(b)=\gamma(\kappa_d b)
 \qquad(b\in \mathcal O_d).
\end{equation}
Then $\mathcal V^\omega=\bigoplus_d\mathcal V_d^\omega$ is an injective homomorphism of graded rings, restricting to the identity on $\mathcal O$, and
\[
 \operatorname{im}(\mathcal V_d^\omega)=\mathcal A_{\mathrm{tr}}[d].
\]
In particular $\mathcal A_{\mathrm{tr}}=\bigoplus_d\mathcal A_{\mathrm{tr}}[d]$ is a graded unital subring and $\mathcal V^\omega:\mathcal A\xrightarrow{\sim}\mathcal A_{\mathrm{tr}}$. Replacing $\omega$ by $c\omega$, $c\in\mathcal O^\times$, replaces this isomorphism by $\operatorname{Ad}(c)\circ\mathcal V^\omega$, where $\operatorname{Ad}(c)(\gamma)=c\gamma c^{-1}$ uses the source and target coefficient actions above.
\end{theorem}

\begin{proof}
Multiplication by $\kappa_d$ is invertible after inverting $p$, so its transpose is injective on the $p$-torsion-free module $\mathcal A[d]$. For $z\in C_{a,b}$, formula \eqref{eq:8.8} and the second identity in \eqref{eq:8.11} give
\[
 \langle\mathcal V_a^\omega\gamma,\mathcal V_b^\omega\eta\rangle(z)
 =\langle\gamma,\eta\rangle((\kappa_a\otimes\kappa_b)z),
\]
hence $\mathcal V_a^\omega(\gamma)\mathcal V_b^\omega(\eta) =\mathcal V_{a+b}^\omega(\gamma\eta)$. Since $\kappa_0=1$, the map is unital and fixes $\mathcal O$. Proposition~\ref{prop:5.2} writes every functional uniquely as $\lambda_d(v-)$, whose image is
\[
 \lambda_d(v\kappa_d-)
 =p^{-d}\operatorname{Tr}(e_d v\kappa_d-)
 =\operatorname{Tr}(\xi_d^{-1}v-).
\]
As $\xi_d$ is a unit, these are exactly the trace functionals. Finally \eqref{eq:8.9} gives $h'_d=t_d(c)s_d(c)^{-1}h_d$ and $\kappa'_d=s_d(c)t_d(c)^{-1}\kappa_d$. Substitution in \eqref{eq:8.15} proves the conjugation formula.
\end{proof}

Notice that for $n\geq2$, after inverting $p$, the algebra $C_{1,1}$ has ranks $1$ and $p+1$ over the loci of subgroups isomorphic to $\mathbf Z/p^2$ and $(\mathbf Z/p)^2$, respectively: these count the subgroups of order $p$ in the kernel of the composite isogeny. Thus it cannot be finite flat over the local ring $\mathcal O_2$, where finite flat modules are free of constant rank.

Nevertheless, after inverting $p$, both $\mathcal O_{a+b}$ and $C_{a,b}$ are finite etale over the source by \eqref{eq:5.7} and base change. A map between finite etale algebras is finite etale, as one sees after an etale cover splitting both into products. Hence the relative trace
\[
 \operatorname{Tr}_{C_{a,b}[1/p]/\mathcal O_{a+b}[1/p]}
\]
is defined. After a base change that splits these finite etale algebras, its value on $1$ at a component indexed by a subgroup $H$ counts the intermediate subgroups of order $p^a$ in $H$.

Using the perfect pairings constructed below, we define the adjoint of $\phi_{a,b}$ over $\mathcal O$. The identity \eqref{eq:8.11} identifies it with this relative trace after inverting $p$, and reduction at the closed point shows that its value on $1$ is a unit.

\begin{theorem}
\label{thm:8.5}
For $a,b\geq0$, the relative trace displayed above extends uniquely to an $\mathcal O_{a+b}$-linear map
\[
 \mathfrak t_{a,b}:C_{a,b}\longrightarrow \mathcal O_{a+b}.
\]
The element $w_{a,b}=\mathfrak t_{a,b}(1)$ is a unit, and $w_{a,b}^{-1}\mathfrak t_{a,b}$ is a $\mathcal O_{a+b}$-linear retraction of $\phi_{a,b}$. For $x\in \mathcal O_a$ and $y\in \mathcal O_b$,
\begin{equation}\label{eq:8.18}
 \boxed{
 \mathcal T_a(x)\mathcal T_b(y)
 =\mathcal T_{a+b}\bigl(\mathfrak t_{a,b}(x\otimes y)\bigr).
 }
\end{equation}
Consequently, with products of submodules meaning finite sums of products,
\begin{equation}\label{eq:8.19}
 \boxed{
 \mathcal A_{\mathrm{tr}}[a]\mathcal A_{\mathrm{tr}}[b]
 =\mathcal A_{\mathrm{tr}}[a+b].
 }
\end{equation}
\end{theorem}

\begin{proof}
For every $d\geq0$, the functional
\[
 \lambda_d^\omega=p^{-d}\operatorname{Tr}_{\mathcal O_d/\mathcal O}(h_d-)
\]
is integral and perfect by Proposition~\ref{prop:5.2}, since $h_d=\xi_de_d$. Now set $d=a+b$. Base change of $\lambda_b^\omega$ along $t_a$ gives a perfect $\mathcal O_a$-linear functional on $C=C_{a,b}$:
\[
 \lambda_{C/\mathcal O_a}^\omega(x\otimes y)
 =x\,t_a(\lambda_b^\omega(y)).
\]
Source linearity of $\lambda_b^\omega$ makes this respect the tensor relation. Composing with $\lambda_a^\omega$ gives
\[
 \lambda_C^\omega(x\otimes y)
 =\lambda_a^\omega\bigl(x\,t_a(\lambda_b^\omega(y))\bigr).
\]
The perfect pairings defined by these functionals give the isomorphisms
\[
 C\xrightarrow{\sim}\operatorname{Hom}_{\mathcal O_a}(C,\mathcal O_a)
 \xrightarrow{\sim}\operatorname{Hom}_{\mathcal O}(C,\mathcal O).
\]
The first sends $z$ to $\lambda_{C/\mathcal O_a}^\omega(z-)$; the second composes with $\lambda_a^\omega$. Their composite proves perfectness of $\lambda_C^\omega$. Trace transitivity and \eqref{eq:8.11} give
\[
 \lambda_C^\omega(z)
 =p^{-d}\operatorname{Tr}_{C/\mathcal O}(\phi_{a,b}(h_d)z).
\]

For each $z\in C$, the expression $b\mapsto\lambda_C^\omega(z\phi_{a,b}(b))$ is source $\mathcal O$-linear. Perfectness of $\lambda_d^\omega$ therefore defines a unique map $\mathfrak t:C\to \mathcal O_d$, adjoint to $\phi_{a,b}$ for these pairings, by
\begin{equation}\label{eq:8.21}
 \lambda_d^\omega(\mathfrak t(z)b)
 =\lambda_C^\omega(z\phi_{a,b}(b))
 \qquad(z\in C,\ b\in \mathcal O_d).
\end{equation}
Applying \eqref{eq:8.21} to $\phi_{a,b}(b_0)z$ and to $z,b_0b$ gives $\mathfrak t(\phi_{a,b}(b_0)z)=b_0\mathfrak t(z)$ by perfectness, proving $\mathcal O_d$-linearity.

After inverting $p$, ordinary trace transitivity identifies the right side of \eqref{eq:8.21} with $\lambda_d^\omega(\operatorname{Tr}_{C[1/p]/\mathcal O_d[1/p]}(z)b)$. Thus $\mathfrak t$ extends this relative trace; $p$-torsion-freeness of $\mathcal O_d$ makes the extension unique. In particular we see that this map is independent of $\omega$.

Put $w=\mathfrak t(1)$. Trace transitivity gives $\operatorname{Tr}_{\mathcal O_d/\mathcal O}(w)=r_a r_b$. Proposition~\ref{prop:3.1} makes $R_d=\mathcal O_d/\mathfrak m\mathcal O_d$ a local Artinian algebra of source $k$-dimension $r_d$ with residue field $k$. In $R_d$, the image of $w$ differs from its scalar residue $\overline w\in k$ by a nilpotent. Nilpotent multiplication has trace zero, so reduction of the source trace gives
\[
 r_d\overline w=r_a r_b\quad\text{in }k.
\]
The rank series \eqref{eq:5.1} gives $r_a\equiv r_b\equiv r_d\equiv1\pmod p$, so $\overline w=1$. Thus $w$ is a unit. By $\mathcal O_d$-linearity, $\mathfrak t(\phi_{a,b}(b))=bw$; this proves the retraction and surjectivity.

The functional paired with $\mathcal T_a(x)\otimes\mathcal T_b(y)$ is $z\mapsto\operatorname{Tr}_{C/\mathcal O}((x\otimes y)z)$, by trace transitivity in $\mathcal O\to\mathcal O_a\to C$. Restricting along $\phi_{a,b}$ and using the relative trace proves \eqref{eq:8.18}. For any $b\in\mathcal O_d$, surjectivity and \eqref{eq:8.3} give $z=\sum_i x_i\otimes y_i$ with $\mathfrak t(z)=b$. Thus $\mathcal T_d(b)=\sum_i\mathcal T_a(x_i)\mathcal T_b(y_i)$, proving the reverse inclusion in \eqref{eq:8.19}; the forward inclusion is already \eqref{eq:8.18}.
\end{proof}

Duality for finite free modules makes this construction commute with source base change. For longer flags, iterating \eqref{eq:8.3} gives finite free source algebras. The two composites of integral traces agree after inverting $p$ and hence over $\mathcal O$ by $p$-torsion-freeness. These identities then base change to arbitrary coefficient rings, including those with $p$-torsion, extending \eqref{eq:8.18} to every finite chain.

\begin{corollary}
Suppose $n\geq2$. Let $\varepsilon_{\mathrm{el}}\in\mathcal O_2[1/p]$ be the idempotent equal to $1$ on the locus of elementary abelian subgroups, isomorphic to $(\mathbf Z/p)^2$, and $0$ on the locus of cyclic subgroups, isomorphic to $\mathbf Z/p^2$. These group types are read after a base change that splits the finite etale groups. Then
\begin{equation}\label{eq:8.22}
 w_{1,1}=1+p\varepsilon_{\mathrm{el}}\in \mathcal O_2^\times,
 \qquad
 \{c\in\mathcal O:c\varepsilon_{\mathrm{el}}\in \mathcal O_2\}=p\mathcal O.
\end{equation}
The source scalars making $1-\varepsilon_{\mathrm{el}}$, the idempotent selecting cyclic subgroups, integral form the same ideal.
\end{corollary}

\begin{proof}
A cyclic group of order $p^2$ has one subgroup of order $p$, whereas $(\mathbf Z/p)^2$ has $p+1$. After inverting $p$, the relative trace of $1$ counts these flags, proving the first identity. Thus $v=p\varepsilon_{\mathrm{el}}$ is integral. Both types occur when $n\geq2$. The idempotent $\varepsilon_{\mathrm{el}}$ cannot lie in the local ring $\mathcal O_2$, so $v\notin p\mathcal O_2$. If $c\varepsilon_{\mathrm{el}}$ is integral, then $cv\in p\mathcal O_2$. The module $\mathcal O_2/p\mathcal O_2$ is free over the domain $\mathcal O/p\mathcal O$; its nonzero element $\overline v$ has zero scalar annihilator. Hence $c\in p\mathcal O$. The reverse inclusion follows from $v\in \mathcal O_2$, and replacing $\varepsilon_{\mathrm{el}}$ by $1-\varepsilon_{\mathrm{el}}$ does not change the ideal of source scalars making the idempotent integral.
\end{proof}

Thus the differing ranks after inverting $p$ that obstruct flatness also produce the integral element $p\varepsilon_{\mathrm{el}}$. At height two, after inverting $p$, the locus of elementary abelian subgroups consists precisely of the subgroup $G[p]$, so \eqref{eq:8.22} recovers Theorem~\ref{thm:9.2} at $d=2,m=1$.

\section{Subgroup congruences from integral trace duality}
\label{sec:9}

For $m\geq0$ and $d\geq nm$, consider the locus of order-$p^d$ subgroups containing $G[p^m]$. Theorem~\ref{thm:9.2} determines which source scalars make its generic projector integral and gives a fibre-product presentation of $\mathcal O_d$ describing the congruence with the complementary locus. Throughout, scalars and duals use the source structures, and ``generic'' means after inverting $p$.

\subsection{The full-torsion locus and its projector}

For $r\geq0$, write
$\operatorname{Sub}_r(G)=\operatorname{Spf}(\mathcal O_r)$ for the
formal scheme classifying finite locally free subgroups of $G$ of
rank $p^r$.

Fix integers $m\geq0$ and $d\geq nm$, and put $\ell=d-nm$. Write $G[p^m]=\ker([p^m]:G\to G)$. For a subgroup $K$ of rank $p^\ell$, its inverse image $p^{-m}K=[p^m]^{-1}(K)$ has rank $p^d$ and contains $G[p^m]$. This construction gives a closed immersion
\[
 p^{-m}:\operatorname{Sub}_{\ell}(G)
 \hookrightarrow\operatorname{Sub}_{d}(G)
\]
whose image is the locus of subgroups containing $G[p^m]$ \cite[Lemma~12.9]{StricklandSubgroups}. On rings write
\[
 \rho_{d,m}:\mathcal O_d\twoheadrightarrow \mathcal O_{\ell},
 \qquad K_m=\ker(\rho_{d,m}).
\]
Inverse image under $[p^m]$ preserves the underlying deformation, so
\[
 \rho_{d,m}\circ s_d=s_\ell.
\]
Thus $\rho_{d,m}$ is a map of source $\mathcal O$-algebras. After inverting $p$, both subgroup algebras are finite etale by \eqref{eq:5.7}, so this closed locus is also open. Its idempotent is
\[
 \varepsilon_m\in \mathcal O_d[1/p],
\]
equal to $1$ on subgroups containing $G[p^m]$ and $0$ on the complement. Multiplication by $\varepsilon_m$ is the generic projector onto this locus.

\begin{lemma}
\label{lem:9.1}
There is a unit $u_m\in \mathcal O_{\ell}^{\times}$ such that
\begin{equation}\label{eq:9.3}
 \rho_{d,m}(e_d)=u_m p^m e_{\ell}.
\end{equation}
\end{lemma}

\begin{proof}
Choose an oriented coordinate $X$ on $G$, and let $f_H(X)$ be the monic polynomial cutting out a finite subgroup divisor $H$. Its constant term is zero because $H$ contains the identity. For $H$ of rank $p^r$, its Euler class is the pullback of $e_r$ under its classifying map; it agrees with $f_H'(0)$ up to a unit \cite[proof of Theorem~9.2]{StricklandSymmetric}, with coefficient extension as in Remark~\ref{rem:5.3}. For the trivial subgroup, $f_H(X)=X$ and $e_0=1$.

Let $K$ be the universal subgroup of order $p^\ell$ over $\mathcal O_\ell$. Strickland's inverse-image construction gives
\[
 f_{p^{-m}K}(X)=v(X)f_K([p^m](X)),
 \qquad v(X)\in \mathcal O_\ell[[X]]^\times
\]
\cite[Definition~12.8]{StricklandSubgroups}. Differentiating at the identity gives
\[
 f_{p^{-m}K}'(0)=v(0)\,p^m f_K'(0).
\]
Here $f_K(0)=0$ eliminates the term involving $v'(0)$, and $[p^m]'(0)=p^m$ since integer multiplication acts by that integer on the tangent line. The factor $v(0)$ and the two Euler comparison factors are units, giving \eqref{eq:9.3}.
\end{proof}

For an ideal $K\subseteq \mathcal O_d$, write
\[
 \operatorname{Ann}_{\mathcal O_d}(K)
 =\{b\in \mathcal O_d:bk=0\text{ for every }k\in K\}.
\]
After inverting $p$, $\operatorname{Ann}(K_m)$ is the selected factor $\varepsilon_m\mathcal O_d[1/p]$; its quotient retains the complement.

\begin{theorem}
\label{thm:9.2}
Let $m\geq0$, $d\geq nm$ and $\ell=d-nm$. Then
\begin{equation}\label{eq:9.5}
 \boxed{
 \{c\in\mathcal O:c\varepsilon_m\in \mathcal O_d\}
 =p^{(n-1)m}\mathcal O.
 }
\end{equation}
Moreover, if
\[
 D_{d,m}=\mathcal O_d/\operatorname{Ann}_{\mathcal O_d}(K_m),
\]
then $D_{d,m}$ is finite free over $\mathcal O$ and there is a canonical fibre-product presentation
\begin{equation}\label{eq:9.6}
 \boxed{
 \mathcal O_d\simeq
 \mathcal O_{\ell}
 \times_{\mathcal O_{\ell}/p^{(n-1)m}\mathcal O_{\ell}}
 D_{d,m}.
 }
\end{equation}
\end{theorem}

To specify the arrows in \eqref{eq:9.6}, put $A_{d,m}=\mathcal O_\ell/p^{(n-1)m}\mathcal O_\ell$. The two maps to the common quotient are
\[
 \begin{aligned}
 \alpha:\mathcal O_\ell&\longrightarrow A_{d,m},
 &x&\longmapsto x\bmod p^{(n-1)m},\\
 \beta:D_{d,m}&\longrightarrow A_{d,m},
 &[b]&\longmapsto \rho_{d,m}(b)\bmod p^{(n-1)m}.
 \end{aligned}
\]
The proof shows that $\beta$ is well-defined and that $b\mapsto(\rho_{d,m}(b),[b])$ identifies $\mathcal O_d$ with the pairs satisfying $\alpha(x)=\beta([b])$.

\begin{proof}
\emph{Step 1: identify the annihilator by duality.} Use the perfect functionals of Proposition~\ref{prop:5.2} at indices $d$ and $\ell$. There is a unique $\eta_m\in \mathcal O_d$ satisfying
\begin{equation}\label{eq:9.8}
 \lambda_d(\eta_m b)=\lambda_{\ell}(\rho_{d,m}(b))
 \qquad(b\in \mathcal O_d).
\end{equation}
Lemma~\ref{lem:5.1}, applied to $K_m$ with finite free quotient $\mathcal O_\ell$, gives the isomorphism
\[
 \operatorname{Ann}_{\mathcal O_d}(K_m)
 \longrightarrow\operatorname{Hom}_{\mathcal O,s_\ell}(\mathcal O_{\ell},\mathcal O),
 \qquad
 a\longmapsto(\rho_{d,m}(b)\mapsto\lambda_d(ab)).
\]
Since every functional on $\mathcal O_\ell$ is $\lambda_\ell(x-)$ and $x$ lifts through $\rho_{d,m}$, equation~\eqref{eq:9.8} identifies its representing element with a multiple of $\eta_m$. Hence
\begin{equation}\label{eq:9.9}
 \operatorname{Ann}_{\mathcal O_d}(K_m)=\eta_m\mathcal O_d.
\end{equation}
The surjection $\rho_{d,m}$ splits as a map of finite free source $\mathcal O$-modules. Its dual injection therefore has direct-summand image, which the perfect pairing identifies with $\operatorname{Ann}(K_m)$.

\smallskip
\noindent\emph{Step 2: compute the generic projector.} Put $q=(n-1)m$. Since the subgroup rings are local (Section~\ref{sec:2}), lift $u_m$ to a unit $\widetilde u_m\in\mathcal O_d^\times$. After inverting $p$, $\rho_{d,m}$ identifies $\varepsilon_m\mathcal O_d[1/p]$ with $\mathcal O_\ell[1/p]$. Trace on this product factor, formula~\eqref{eq:5.4} and Lemma~\ref{lem:9.1} give
\[
 \begin{aligned}
 \lambda_d(\widetilde u_m^{-1}p^q\varepsilon_m b)
 &=p^{q-d}\operatorname{Tr}_{\mathcal O_\ell/\mathcal O}
       (u_m^{-1}\rho_{d,m}(e_d)\rho_{d,m}(b))\\
 &=p^{-\ell}\operatorname{Tr}_{\mathcal O_\ell/\mathcal O}
       (e_\ell\rho_{d,m}(b))
 =\lambda_\ell(\rho_{d,m}(b)),
 \end{aligned}
\]
where traces and functionals are extended to $\mathcal O[1/p]$ and $q-d+m=-\ell$. Perfectness persists under this base change, so comparison with \eqref{eq:9.8} gives
\begin{equation}\label{eq:9.10}
 \eta_m=\widetilde u_m^{-1}p^{(n-1)m}\varepsilon_m
 \qquad\text{in }\mathcal O_d[1/p].
\end{equation}
Its integral left side proves $p^q\varepsilon_m\in\mathcal O_d$.

\smallskip
\noindent\emph{Step 3: prove that the scalar denominator is exact.} If $c\varepsilon_m\in\mathcal O_d$ for $c\in\mathcal O$, it annihilates $K_m$: this holds after inverting $p$ and descends by $p$-torsion-freeness. Thus it lies in \eqref{eq:9.9}. Applying $\rho_{d,m}$ and using \eqref{eq:9.10} gives $c\in p^q\mathcal O_\ell$. The algebra $\mathcal O_\ell$ is finite free of positive rank, hence faithfully flat; ideal contraction gives $c\in p^q\mathcal O$, proving \eqref{eq:9.5}.

\smallskip
\noindent\emph{Step 4: recover the integral algebra by gluing.} The intersection $K_m\cap\operatorname{Ann}(K_m)$ vanishes after inverting $p$, hence vanishes integrally by $p$-torsion-freeness. Thus $b\mapsto(\rho_{d,m}(b),[b])$ is injective. Equations~\eqref{eq:9.9} and~\eqref{eq:9.10} give
\[
 \rho_{d,m}(\operatorname{Ann}(K_m))=p^q\mathcal O_\ell,
\]
making $\beta$ well-defined. For a pair $(x,[b])$ with matching reductions, choose $a\in\operatorname{Ann}(K_m)$ with $\rho_{d,m}(a)=x-\rho_{d,m}(b)$; then $b+a$ lifts the pair. This proves \eqref{eq:9.6}. The direct-summand statement in Step~1 makes $D_{d,m}$ finite projective, hence free, over local $\mathcal O$.
\end{proof}

Since $\alpha$ is surjective, the fibre product gives the exact sequence of source $\mathcal O$-modules
\[
 0\longrightarrow \mathcal O_d
 \xrightarrow{\,b\mapsto(\rho_{d,m}(b),[b])\,}
 \mathcal O_\ell\oplus D_{d,m}
 \xrightarrow{\,\alpha-\beta\,} A_{d,m}
 \longrightarrow0.
\]
The cokernel $A_{d,m}$ is the \emph{congruence module}: it measures the failure of independent functions on the two pieces to extend to an integral function on the whole subgroup scheme.

For $m=0$, $\rho_{d,0}$ is the identity, $\varepsilon_0=1$ and $D_{d,0}=0$. At height one the unique subgroup of each generic order contains $G[p^m]$ when $m\leq d$, so $\varepsilon_m=1$ and $K_m[1/p]=0$. Torsion-freeness gives $K_m=0$ and $D_{d,m}=0$. In both cases the common quotient is zero and the fibre product reduces to $\mathcal O_\ell$; the complementary module is free of rank zero.

\subsection{Individual subgroup types at height two}

\begin{corollary}
\label{cor:9.3}
Assume $n=2$. By \eqref{eq:5.7} and the invariant-factor classification of finite abelian groups, each geometric generic subgroup of order $p^d$ has a unique type (abstract group isomorphism class)
\[
 H_a\simeq\mathbf Z/p^{d-a}\oplus\mathbf Z/p^a,
 \qquad 0\leq a\leq L:=\lfloor d/2\rfloor.
\]
A zero exponent denotes the trivial factor. Let $\pi_a\in \mathcal O_d[1/p]$ be the idempotent selecting subgroups of type $H_a$. Then
\[
 \boxed{
 \min\{s\in\mathbf Z_{\geq0}:p^s\pi_a\in \mathcal O_d\}
 =
 \begin{cases}
 a+1,&a<L,\\
 a,&a=L.
 \end{cases}
 }
\]
In particular, the cyclic projector has exact denominator $p$ whenever $d\geq2$.
\end{corollary}

\begin{proof}
For $0\leq a\leq L$, the earlier projector $\varepsilon_a$ selects subgroups containing $G[p^a]$, with $\varepsilon_0=1$. Since $d\geq2a$, Theorem~\ref{thm:9.2} gives
\[
 \min\{s\in\mathbf Z_{\geq0}:p^s\varepsilon_a\in \mathcal O_d\}=a.
\]
Containment of $G[p^a]$ is equivalent to the smaller invariant exponent being at least $a$: precisely then $H[p^a]\subseteq G[p^a]$ has the same order $p^{2a}$ and equals it. Thus
\[
 \pi_a=\varepsilon_a-\varepsilon_{a+1}\quad(a<L),
 \qquad
 \pi_L=\varepsilon_L.
\]
For $a<L$, also $a+1\leq L$, so multiplication by $p^{a+1}$ makes both terms integral. If $p^a\pi_a$ were integral, then $p^a\varepsilon_{a+1}=p^a\varepsilon_a-p^a\pi_a$ would be integral, contradicting its exact denominator $p^{a+1}$. The final projector is $\varepsilon_L$.
\end{proof}

For example, for $d=2$ the cyclic and elementary abelian projectors both have denominator $p$. The elementary type is the full-torsion locus $G[p]$, so Theorem~\ref{thm:9.2}, with $m=1$, gives the common quotient $\mathcal O_0/p\mathcal O_0=\mathcal O/p$. For $d=4$, the types with invariant exponents $(4,0)$, $(3,1)$ and $(2,2)$ have denominator exponents $1,2,2$, respectively. For $d=0$ and $d=1$ there is a single type, whose projector is $1$ and whose denominator exponent is zero.

\begin{remark}
Theorem~\ref{thm:9.2} determines the congruence joining the full-torsion locus to its complement. Corollary~\ref{cor:9.3} determines individual projector denominators in height two. A presentation of $\mathcal O_d$ specifying the simultaneous integral
congruences between all generic subgroup types in height $n\geq3$
is not obtained here.
\end{remark}

\section*{Statement on the use of generative AI}

Generative AI tools were used during the development of this manuscript to help locate references, assist with reviewing
the proofs of Theorems~\ref{thm:5.10} and~\ref{thm:7.5} and Lemma~\ref{lem:5.5}, and identify unclear or potentially misleading passages in the exposition.

The author takes full responsibility for the mathematical content of the paper, including the correctness of its statements and proofs and the accuracy of its references.


\begin{thebibliography}{99}
\raggedright

\bibitem{AdepuChoudharyKim2026} A.~Adepu, N.~Choudhary, and G.~Kim, \emph{Finite Quotient Towers and Geometric Witt Filtrations}, preprint, 19 October 2026; \url{https://github.com/Quasar0147/witt-filtrations/blob/b3ba5210f302e2f3961117498895b0154d16e6a8/preprint.pdf}.

\bibitem{AHSSigma} M.~Ando, M.~J.~Hopkins, and N.~P.~Strickland, \emph{The sigma orientation is an $H_\infty$ map}, Amer. J. Math. \textbf{126} (2004), no.~2, 247--334; arXiv:math/0204053v2.

\bibitem{BehrensConstruction} M.~Behrens, \emph{The construction of tmf}, in \emph{Topological Modular Forms} (C.~L.~Douglas, J.~Francis, A.~G.~Henriques, and M.~A.~Hill, eds.), Math. Surveys Monogr. \textbf{201}, Amer. Math. Soc., Providence, RI, 2014, Chapter~12, 131--188; \url{https://www.math.ucla.edu/~mikehill/Research/surv-douglas2-201.pdf}.

\bibitem{BSY} R.~Burklund, T.~M.~Schlank, and A.~Yuan, \emph{The chromatic Nullstellensatz}, arXiv:2207.09929, 2022; to appear in Ann. of Math. Numbering follows arXiv:2207.09929v1.

\bibitem{LangAlgebra} S.~Lang, \emph{Algebra}, revised third edition, Graduate Texts in Mathematics \textbf{211}, Springer, New York, 2002; \url{https://doi.org/10.1007/978-1-4613-0041-0}.

\bibitem{LazardFormalGroups} M.~Lazard, \emph{Sur les groupes de Lie formels \`a un param\`etre}, Bull. Soc. Math. France \textbf{83} (1955), 251--274; \url{https://doi.org/10.24033/bsmf.1462}.

\bibitem{LubinTateFormalModuli} J.~Lubin and J.~Tate, \emph{Formal moduli for one-parameter formal Lie groups}, Bull. Soc. Math. France \textbf{94} (1966), 49--59; \url{https://doi.org/10.24033/bsmf.1633}.

\bibitem{NelsonWreath} P.~Nelson, \emph{On the Morava $E$-theory of wreath products of symmetric groups}, arXiv:1803.05020v1, 2018.

\bibitem{RezkCongruence} C.~Rezk, \emph{The congruence criterion for power operations in Morava $E$-theory}, Homology Homotopy Appl. \textbf{11} (2009), no.~2, 327--379; arXiv:0902.2499v2.

\bibitem{RezkHeightTwo} C.~Rezk, \emph{Power operations for Morava $E$-theory of height $2$ at the prime $2$}, \url{https://arxiv.org/abs/0812.1320v1}, 2008.

\bibitem{RezkKoszul} C.~Rezk, \emph{Rings of power operations for Morava $E$-theories are Koszul}, arXiv:1204.4831v2; author version dated 10 July 2017, \url{https://rezk.web.illinois.edu/rezk-koszul-morava-e-new-version.pdf}.

\bibitem{RezkIsogenies} C.~Rezk, \emph{Isogenies, power operations, and homotopy theory}, in \emph{Proceedings of the International Congress of Mathematicians---Seoul 2014}, Vol.~II, Kyung Moon Sa, Seoul, 2014, 1125--1145.

\bibitem{RezkLogarithm} C.~Rezk, \emph{The units of a ring spectrum and a logarithmic cohomology operation}, J. Amer. Math. Soc. \textbf{19} (2006), no.~4, 969--1014; arXiv:math/0407022v2.

\bibitem{RezkWitt} C.~Rezk, \emph{The Witt filtration of Lubin--Tate deformation rings}, arXiv:2603.12490v1, 2026.

\bibitem{StapletonTranschromatic} N.~Stapleton, \emph{Transchromatic generalized character maps}, Algebr. Geom. Topol. \textbf{13} (2013), no.~1, 171--203; arXiv:1110.3346v2.

\bibitem{StapletonSubgroups} N.~Stapleton, \emph{Subgroups of $p$-divisible groups and centralizers in symmetric groups}, Trans. Amer. Math. Soc. \textbf{367} (2015), no.~5, 3733--3757; \url{https://doi.org/10.1090/S0002-9947-2014-06344-7}. Numbering follows arXiv:1305.6573v1.

\bibitem{StricklandDuality} N.~P.~Strickland, \emph{$K(n)$-local duality for finite groups and groupoids}, Topology \textbf{39} (2000), no.~4, 733--772; arXiv:math/0011109v1.

\bibitem{StricklandSubgroups} N.~P.~Strickland, \emph{Finite subgroups of formal groups}, J. Pure Appl. Algebra \textbf{121} (1997), no.~2, 161--208. Theorem and section numbering follows the author's online version, \url{https://strickland1.org/research/subgp.pdf}.

\bibitem{StricklandSymmetric} N.~P.~Strickland, \emph{Morava $E$-theory of symmetric groups}, Topology \textbf{37} (1998), no.~4, 757--779; arXiv:math/9801125v1.

\bibitem{StacksFlatness} The Stacks Project Authors, \emph{The Stacks Project}, Lemma~10.99.7, Local criterion for flatness, \url{https://stacks.math.columbia.edu/tag/00MK}.

\bibitem{StacksDescent} The Stacks Project Authors, \emph{The Stacks Project}, Lemmas~35.23.32 and~35.37.1, Descent of finite locally free and affine morphisms, \url{https://stacks.math.columbia.edu/tag/02VO}, \url{https://stacks.math.columbia.edu/tag/0245}.

\bibitem{ZhuHecke} Y.~Zhu, \emph{The Hecke algebra action and the Rezk logarithm on Morava $E$-theory of height $2$}, Trans. Amer. Math. Soc. \textbf{373} (2020), no.~5, 3733--3764; \url{https://yifeizhu.github.io/ho.pdf}.

\end{thebibliography}
\end{document}